\documentclass[11pt]{article}
\usepackage[margin=1in]{geometry}
\usepackage{amsmath,amssymb,amsthm,mathtools}
\usepackage{enumitem}
\usepackage{microtype}
\usepackage{xcolor}
\usepackage[numbers,sort&compress]{natbib}
\usepackage[
  colorlinks=true,
  linkcolor=blue!55!black,
  citecolor=blue!55!black,
  urlcolor=blue!55!black,
  pdftitle={When Weak Submodularity Meets Matroids: How Bad Can Standard Greedy Be?},
  pdfauthor={Qixin Zhang},
  pdfsubject={Approximation guarantees for weakly submodular maximization under matroid constraints},
  pdfkeywords={weak submodularity, matroid constraint, greedy algorithm, approximation ratio, combinatorial optimization}
]{hyperref}

\newtheorem{theorem}{Theorem}[section]
\newtheorem{lemma}[theorem]{Lemma}
\newtheorem{proposition}[theorem]{Proposition}

\newcommand{\R}{\mathbb{R}}
\newcommand{\cI}{\mathcal{I}}
\newcommand{\cM}{\mathcal{M}}
\newcommand{\OPT}{\mathrm{OPT}}
\newcommand{\gr}{\mathrm{gr}}
\newcommand{\eps}{\varepsilon}
\DeclareMathOperator*{\argmax}{arg\,max}

\makeatletter
\newcommand{\appendixcontents}{\@starttoc{apc}}
\newcommand*\l@appsubsection{\@dottedtocline{2}{1.5em}{3.6em}}
\makeatother
\newcommand{\addappendixsection}[1]{%
  \addcontentsline{apc}{section}{\protect\numberline{\thesection}#1}}
\newcommand{\addappendixsubsection}[1]{%
  \addcontentsline{apc}{appsubsection}{\protect\numberline{\thesubsection}#1}}

\title{\bfseries When Weak Submodularity Meets Matroids: How Does Standard Greedy Behave?}
\author{
Qixin Zhang\textsuperscript{1},\quad
Shi Fu\textsuperscript{1},\quad
Shengminjie Chen\textsuperscript{2},\quad
Chao Qian\textsuperscript{3},\quad
Dacheng Tao\textsuperscript{1}
\\[0.8em]
\small
\textsuperscript{1}College of Computing and Data Science,
Nanyang Technological University, Singapore
\\
\small
\textsuperscript{2}Institute of Computing Technology,
Chinese Academy of Sciences, Beijing, China
\\
\small
\textsuperscript{3}School of Artificial Intelligence,
Nanjing University, Nanjing, China
\\[0.6em]
\small
\href{mailto:qixin.zhang2026@gmail.com}
{\texttt{qixin.zhang2026@gmail.com}}
\quad
\href{mailto:shi011@e.ntu.edu.sg}
{\texttt{shi011@e.ntu.edu.sg}}
\\
\small
\href{mailto:csmj@ict.ac.cn}
{\texttt{csmj@ict.ac.cn}}
\quad
\href{mailto:qianc@lamda.nju.edu.cn}
{\texttt{qianc@lamda.nju.edu.cn}}
\quad
\href{mailto:dacheng.tao@ntu.edu.sg}
{\texttt{dacheng.tao@ntu.edu.sg}}
}

\date{}

\begin{document}
\maketitle

\begin{abstract}
Weakly submodular maximization finds numerous applications in machine learning, data mining and statistics.  In particular, under a cardinality constraint, the standard greedy algorithm
is known to attain the tight $(1-e^{-\gamma})$-approximation guarantee,
where $\gamma$ denotes the submodularity ratio. However, when considering  general matroid constraints, \emph{how well standard greedy can perform} has remained a
longstanding open question. In this paper, we resolve this question by studying the worst-case approximation guarantee $R_{\mathrm{gr}}(\gamma,r)$ of standard forward greedy, where $r$ is the rank of matroid. Specifically, our first contribution is the following refined lower bound
\begin{equation*}
R_{\gr}(\gamma,r)\ge \max\left\{
\frac{\gamma}{r},
\frac{\gamma^2}
{B_\gamma+\sqrt{\gamma(r-1)}},
\frac{\gamma^{r+1}}
{2-\gamma+\gamma^{r+1}},
\frac{\gamma^2}
{14\left(B_\gamma+\sqrt{\gamma(1-\gamma)(r-1)}\right)}
\right\},
\end{equation*} where $B_\gamma\triangleq2-\gamma+\gamma^2$. This lower bound not only
strictly improves upon the previous guarantee of
$0.4\gamma^2/(1+\sqrt{\gamma r})$ of
\cite{gatmiry2018nonsubmodular}, but also can recover the classical
$1/2$-approximation guarantee at $\gamma=1$. Subsequently, a corresponding upper bound is established for $R_{\mathrm{gr}}(\gamma,r)$, namely,
\begin{equation}\label{eq_upper_abstract}
R_{\gr}(\gamma,r)
\le\min\left\{\frac{1}{2},20\max\left\{\frac{\gamma}{r},
\frac{\gamma^2}{B_\gamma+
\sqrt{\gamma(1-\gamma)(r-1)}}\right\}\right\}.
\end{equation} For every fixed $\gamma\in(0,1)$, this Eq.\eqref{eq_upper_abstract} effectively sharpens the recent
$O_\gamma(1/\log(n-1))$ upper bound of \cite{feldman2026maximizing} to $O_\gamma(r^{-1/2})$, where $n$ is the size of ground set. Putting all these results together, we eventually can obtain
\[R_{\gr}(\gamma,r)
=
\Theta\!\left(
\max\left\{\frac{\gamma}{r},
\frac{\gamma^2}
{1+\sqrt{\gamma(1-\gamma)r}}\right\}
\right).\]To the best of our knowledge, this is the first asymptotically tight
characterization of the worst-case performance of standard greedy
for matroid-constrained weakly submodular maximization

\end{abstract}

\section{Introduction}

Submodularity arises naturally in a wide range of applications, including influence
maximization~\cite{kempe2003influence,borgs2014maximizing,tang2018online,zhang2026one}, feature
selection~\cite{qian2015subset,qian2017subset,elenberg2018restricted}, sensor placement~\cite{krause2008sensor,hashemi2019submodular}, and data
summarization~\cite{lin2011summarization,mirzasoleiman2013distributed,badanidiyuru2014streaming}. A standard approach for maximizing submodular objectives is forward greedy: starting from the empty set, the
algorithm repeatedly adds a feasible element with the largest marginal gain.
Furthermore, classical theoretical analysis also gives that, for monotone submodular functions,
forward greedy can achieve a tight approximation ratio of $1-1/e$ under a
cardinality constraint~\cite{nemhauser1978analysis}.

However, exact submodularity may be overly restrictive for many practical applications. Indeed,
many real-world scenarios give rise to utility functions that are close to submodular. Important examples include 
variable selection for regression~\cite{das2011submodular,das2018approximate,khanna2017scalable}, neural network pruning~\cite{elenberg2018restricted} and sparse optimal transport~\cite{manupriya2024submodular}. One of the most widely used notions for characterizing such  close-to-submodular objectives is
$\gamma$-weak submodularity~\cite{das2011submodular,das2018approximate,qian2018approximation}. More specifically, a monotone set function $f:2^U\to\R_+$ is said to be $\gamma$-weakly
submodular if and only if, for every $A\subseteq C\subseteq U$,
\begin{equation*}
\sum_{e\in C\setminus A}\left(f(A\cup\{e\})-f(A)\right)
\ge \gamma\bigl(f(C)-f(A)\bigr).
\end{equation*}This inequality quantifies the extent to which the joint gain $f(C)-f(A)$ can exceed the sum of
these individual marginal gains obtained
by adding each element of $C\setminus A$ to
$A$. In particular, every monotone submodular
function satisfies the above inequality with $\gamma=1$. More importantly, under a cardinality constraint, Das and Kempe~\cite{das2011submodular,das2018approximate}
proved that standard forward greedy achieves a
$(1-e^{-\gamma})$-approximation for monotone
$\gamma$-weakly submodular maximization. This guarantee was later shown to be
tight in the value-oracle model by \cite{harshaw2019submodular}.

Beyond the cardinality constraints considered above, many of real-world applications also require us to obey additional structural restrictions, such as diversity
and fairness, which can be naturally modeled by matroid constraints. For normalized monotone submodular functions, Fisher et
al.~\cite{fisher1978analysis} showed that standard forward greedy can attain a $1/2$-approximation guarantee under a matroid. As for weakly submodular objectives,
however, the approximation behavior of standard forward greedy remains far
less understood. On the one hand, Gatmiry and Gomez
Rodriguez~\cite{gatmiry2018nonsubmodular} established an approximation
guarantee of $ \frac{0.4\gamma^2}{1+\sqrt{\gamma r}}$. On the other hand, recent work~\cite{feldman2026maximizing} constructed a family of bad
instances showing that, for every fixed $\gamma\in(0,1)$, the approximation
ratio of standard forward greedy can be as small as $O(\frac{\gamma}{(1-\gamma)^{2}\log(n-1)})$, where $n$ is the size of ground set. Consequently, there exists a substantial gap between the
known lower and upper bounds for standard forward greedy on matroid-constrained weakly submodular maximization. In the remainder of this paper, we close this gap up to a constant factor. 

Specifically, letting $R_{\mathrm{gr}}(\gamma,r)$ denote the worst-case
approximation ratio of standard forward greedy over all normalized monotone
$\gamma$-weakly submodular functions and all matroids of rank $r$, we
establish the following matching lower and upper bounds.
\begin{theorem}
\label{thm:complete-bounds}
For every $\gamma\in(0,1]$ and every integer $r\ge 2$, let $B_\gamma\triangleq 2-\gamma+\gamma^2$. Then, the worst-case approximation ratio of standard forward greedy satisfies
\begin{equation*}
R_{\gr}(\gamma,r)
\ge
\max\left\{
    \frac{\gamma}{r},
    \frac{\gamma^2}
    {B_\gamma+\sqrt{\gamma(r-1)}},
    \frac{\gamma^{r+1}}
    {2-\gamma+\gamma^{r+1}},
    \frac{\gamma^2}
    {14\left(B_\gamma+\sqrt{\gamma(1-\gamma)(r-1)}\right)}
\right\},
\end{equation*}
and
\begin{equation*}
R_{\gr}(\gamma,r)
\le
\min\left\{
    \frac{1}{2},
    20\max\left\{
        \frac{\gamma}{r},
        \frac{\gamma^2}
        {B_\gamma+\sqrt{\gamma(1-\gamma)(r-1)}}
    \right\}
\right\}.
\end{equation*}
Consequently,
\begin{equation}
\label{eq:complete-joint-order}
R_{\gr}(\gamma,r)
=
\Theta\!\left(
    \max\left\{
        \frac{\gamma}{r},
        \frac{\gamma^2}
        {1+\sqrt{\gamma(1-\gamma)r}}
    \right\}
\right),
\end{equation}
where the symbol $\Theta(\cdot)$ denotes equality up to multiplicative
constant factors.
\end{theorem}

\subsection{Related Work}
The notion of $\gamma$-weak submodularity was introduced by Das and
Kempe~\cite{das2011submodular,das2018approximate}, who proved that standard
forward greedy achieves a $(1-e^{-\gamma})$-approximation for maximizing a
monotone $\gamma$-weakly submodular function subject to a cardinality
constraint. Bian et al.~\cite{bian2017guarantees} refined this analysis by
incorporating generalized curvature, while Harshaw et
al.~\cite{harshaw2019submodular} showed that the
$(1-e^{-\gamma})$ dependence is tight in the value-oracle model. Furthermore, Qian et
al.~\cite{qian2019maximizing} further showed that the multi-objective
evolutionary algorithm also can attain a
$(1-e^{-\gamma})$-approximation. Disser and Weckbecker~\cite{disser2024unified} subsequently introduced a broader class
of $\gamma$-$\alpha$-augmentable functions,for which standard
forward greedy can obtain tight approximation guarantees under cardinality constraints.

Moving beyond cardinality constraints, Chen et
al.~\cite{chen2018weakly} proved that Residual Random Greedy achieves the
rank-independent approximation factor $ (1+1/\gamma)^{-2}$. To improve this result, Thiery and
Ward~\cite{thiery2022two} then introduced the upper submodularity ratio $\beta$
and developed a distorted local-search algorithm for
$(\gamma,\beta)$-weakly submodular maximization over matroids. Their algorithm achieves an  approximation factor of  $\gamma^2(1-e^{-\phi})/\phi$, where $\phi\triangleq\gamma^2+\beta(1-\gamma)$. Particularly when considering the partition matroids, \cite{zhang2025multinoulli} introduced a novel Multinoulli extension and its associated Multinoulli-SCG algorithm,
which cam preserve the same approximation guarantee of  \cite{thiery2022two} while requiring fewer function
evaluations. Recently, \cite{fu2026poisson}  proposed a
maximum-gain Poisson exchange algorithm that achieves a better approximation ratio
arbitrarily close to $1-
    (\gamma/(2-\gamma))^{\gamma^2/2(1-\gamma)
    }$ for matroid-constrained weakly submodular maximization.  Furthermore, a parallel line of research~\cite{gatmiry2018nonsubmodular,qian2018multiset,gong2019parametric,
lu2022maximizing} explores a stronger weakly DR-submodular maximization problems. 

\section{Preliminaries}

\paragraph{Notation.} Let $U$ be a finite ground set and let  $f:2^U\to\R_+$ be a set function.  For any subset
$S\subseteq U$ and $e\in U$, we write $S+e\triangleq S\cup\{e\}$ and define the
marginal contribution of $e$ to $S$ by $f(e\mid S)\triangleq f(S+e)-f(S)$.  A set
function is \emph{normalized} if $f(\varnothing)=0$, and it is \emph{monotone}
if $f(A)\le f(C)$ for all $A\subseteq C\subseteq U$.  Unless stated otherwise,
we suppose all set functions in this paper are \emph{normalized} and \emph{monotone}.

\paragraph{Weak Submodularity.}For $\gamma\in(0,1]$, a set function $f:2^U\to\R_+$ is \emph{$\gamma$-weakly
submodular} if, for every $A\subseteq C\subseteq U$, the following inequality holds:
\begin{equation}\label{eq:lower-ws}
\sum_{e\in C\setminus A} f(e\mid A)
\ge \gamma\bigl(f(C)-f(A)\bigr).
\end{equation}

\paragraph{Matroids.}A \emph{matroid} is a pair $\cM\triangleq(U,\cI)$, where $U$ is a finite ground set and $\cI\subseteq 2^U$ is a collection of sets satisfying the following
axioms: (i)~$\varnothing\in\cI$; (ii)~if $C\in\cI$ and $A\subseteq C$, then
$A\in\cI$; and (iii)~if $A,C\in\cI$ and $|A|<|C|$, then there exists
$e\in C\setminus A$ such that $A+e\in\cI$.  The members of $\cI$ are called
\emph{independent sets}, and its maximal members are called
\emph{bases}.  Every base has the same cardinality, called the \emph{rank} of
$\cM$ and denoted by $r$. 

\paragraph{Standard Forward Greedy Method.}Starting from $S_0=\varnothing$, at each iteration $t=1,\ldots,r$, the
standard forward greedy algorithm chooses
\begin{equation}\label{eq:greedy-rule}
s_t\in\argmax\{f(e\mid S_{t-1}):
e\notin S_{t-1},\ S_{t-1}+e\in\cI\},
\qquad S_t=S_{t-1}+s_t.
\end{equation}
After $r$ iterations, the algorithm returns the base $S_r$.  The goal of this paper is to
determine the worst-case approximation guarantee of the greedy set $S_r$ under weak submodularity. More specifically, for a fixed $\gamma$-weakly
submodular function $f$,  we  define 
\begin{equation}\label{eq:R-definition}
R_{\gr}(\gamma,r)\triangleq\inf_{(f,\cM)}\frac{f(S_r)}{\max_{T\in\cI}f(T)},
\end{equation}
where the infimum is over all rank-$r$ matroid $\cM$ and all possible objectives satisfying Eq.\eqref{eq:greedy-rule}.  Note that we only consider
instances with positive optimum.

\section{Upper bounds for
\texorpdfstring{$R_{\gr}(\gamma,r)$}{Rgr(gamma,r)}}
\label{sec:upper-bounds}

To characterize upper bounds on $R_{\gr}(\gamma,r)$, we consider a structured family of matroid-constrained weakly submodular maximization instances. Fix positive integers $k$ and $N$, and let the ground set be the disjoint union of three types of elements:
\[
U=G\mathbin{\uplus}B\mathbin{\uplus}E,
\qquad
G=\{g_1,\ldots,g_k\},\quad
B=\{b_1,\ldots,b_k\},\quad
E=\{e_1,\ldots,e_N\}.
\]
Here, $G$ denotes the set of \emph{good} elements, $B$ the set of \emph{bad} elements, and $E$ an \emph{auxiliary} set.

Moreover, we impose a partition matroid consisting of the $k$ bad--good blocks
$\{b_i,g_i\}$ and the $N$ singleton auxiliary blocks $\{e_j\}$, each with capacity
one.  Its rank is therefore $r\triangleq k+N$.

The high-level idea is to design an objective function that forces the greedy algorithm to select all bad elements $\{ b_1,\ldots,b_k\}$. Once a bad element $b_i$ is selected, the partition constraint immediately makes its paired good element $g_i$ infeasible, thereby preventing it from being selected at any later stage. At the same time, the objective is calibrated so that the good elements collectively yield a substantially larger value in an optimal solution. This gap between the greedy solution and the optimal solution forms the key mechanism behind our upper-bound construction.

To achieve this goal, we consider a simplified set objective function that depends only on the numbers of selected elements of each type.
For a set $S\subseteq U$,  we write
\[
b(S)=|S\cap B|,\qquad q(S)=|S\cap G|,
\qquad s(S)=|S\cap E|.
\]
Given a weight $w>0$ and a nonnegative count core
$H:\mathbb Z_+^2\to\R_+$, we define
\begin{equation}\label{eq:common-template}
f_{w,H}(S)
=w\cdot b(S)+H(q(S),s(S)).
\end{equation}
We next design $H$ so that $f_{w,H}$ is $\gamma$-weakly submodular while
the standard forward greedy method has a small approximation ratio.

\subsection{First core: a variant of the harmonic mean}

Recall that the harmonic mean of two positive numbers $x$ and $s$ is
\[
\phi(x,s)=\frac{2xs}{x+s}=\frac{2}{x^{-1}+s^{-1}}.
\]With the help of $\phi(x,s)$, we next define a continuous core $H(q,s)$, that is, 
\begin{equation}\label{eq:H-def}
H(q,s)=\frac12\phi\bigl(a(q+\lambda),s\bigr)
=\frac{a(q+\lambda)s}{a(q+\lambda)+s},
\qquad q,s\ge0,
\end{equation} where $a\ge1$ and
$\lambda>0$ are two parameters to be specified later.
Here $q$ and $s$ will later be instantiated as the numbers of good and
auxiliary elements, respectively. That is,  $\widehat H(S)=H\bigl(|S\cap G|,|S\cap E|\bigr) $ for any $S\subseteq U$. We next state the properties of both core $H$ and $\widehat H$. Their proofs are presented in Appendix~\ref{app:upper-proofs}.

\begin{lemma}[Continuous properties of the harmonic core]
\label{lem:H-continuous}
The function $H$ in Eq.\eqref{eq:H-def} is nonnegative, coordinate-wise nondecreasing, and concave on $\R_{\ge0}^2$.  If
$x\triangleq a(q+\lambda)$ and $D\triangleq x+s$, then
\begin{align}
\partial_q H(q,s)&=\frac{as^2}{D^2},&
H(q+1,s)-H(q,s)&=\frac{D}{D+a}\partial_q H(q,s),
\label{eq:H-forward-q}\\
\partial_s H(q,s)&=\frac{x^2}{D^2},&
H(q,s+1)-H(q,s)&=\frac{D}{D+1}\partial_s H(q,s).
\label{eq:H-forward-s}
\end{align}
\end{lemma}Importantly, with the result of Lemma~\ref{lem:H-continuous}, we can prove that \begin{lemma}
\label{lem:H-ratio}
The set function $\widehat H(S)\triangleq H\bigl(|S\cap G|,|S\cap E|\bigr)$  is monotone $\gamma_\lambda$-weakly submodular, where $\gamma_{\lambda}\triangleq\frac{\lambda}{1+\lambda}$.
\end{lemma}

\subsubsection{How the harmonic core misleads standard forward greedy?}\label{sec:how_the_harmonic}

Throughout this subsubsection, we assume $N\ge k$, set $a=N/k\ge1$, and choose $\lambda\triangleq\frac{\gamma}{1-\gamma}$. Then $\gamma_\lambda=\lambda/(1+\lambda)=\gamma$, so
Lemma~\ref{lem:H-ratio} implies that $\widehat H$ is normalized, monotone, and
$\gamma$-weakly submodular.  Generally speaking, adding the nonnegative modular function
will preserves these properties.  Hence the objective
\begin{equation}\label{eq:offset-objective}
f(S)=|S\cap B|+
H\bigl(|S\cap G|,|S\cap E|\bigr).
\end{equation}
is normalized, monotone, and $\gamma$-weakly submodular.

Whenever $S$ contains neither good nor auxiliary elements,
equivalently $S\subseteq B$, the relevant marginal values are
\[
f(b_i\mid S)=1,\qquad f(g_i\mid S)=0,\qquad
f(e\mid S)=\frac{a\lambda}{a\lambda+1}<1,
\] where $b_i\notin S$. Thus the greedy method strictly chooses all $k$ bad elements first, which will block all good elements. It then
takes all auxiliary elements and returns $S_r\triangleq B\cup E$.  The competing base
$G\cup E$ is feasible. According to  Eq.~\eqref{eq:offset-objective} and $a=N/k$, we can have
\begin{align}
f(S_r)\triangleq k+\frac{\lambda N}{k+\lambda}, \qquad
f(G\cup E)\triangleq\frac{N(k+\lambda)}{2k+\lambda}.
\end{align}
Since $\OPT\triangleq\max_{T\in\cI}f(T)\ge f(G\cup E)$, the worst-case approximation ratio $R_{\gr}(\gamma,r)$ satisfies
\begin{equation}\label{eq:offset-ratio}
\begin{aligned}
R_{\gr}(\gamma,r)\le\frac{f(S_r)}{\OPT}
&\le\frac{2k+\lambda}{k+\lambda}
\left(\frac{k}{N}+\frac{\lambda}{k+\lambda}\right)\le2\left(\frac{k}{N}+\frac{\lambda}{k}\right),
\end{aligned}
\end{equation} where the final inequality comes from
$(2k+\lambda)/(k+\lambda)\le2$ and
$\lambda/(k+\lambda)\le\lambda/k$.  As a result, we naturally have the following
upper bound for $R_{\gr}(\gamma,r)$, that is to say,

\begin{proposition}\label{prop:harmonic-upper}
For every $0<\gamma<1$ and $r\ge2$, let
$\lambda\triangleq\gamma/(1-\gamma)$.  Then, we have 
\begin{equation}\label{eq:harmonic-upper}
R_{\gr}(\gamma,r)
\le
\min_{1\le k\le\lfloor r/2\rfloor}
\frac{2k+\lambda}{k+\lambda}
\left(\frac{k}{r-k}+\frac{\lambda}{k+\lambda}\right).
\end{equation}
Consequently, for every fixed $\gamma\in(0,1)$,
\begin{equation}\label{add_eq:harmonic-fixed-rate}
  R_{\gr}(\gamma,r)
\le\frac1{\sqrt r}
\max\left\{
8\sqrt{\frac{\gamma}{1-\gamma}},
2\sqrt{\frac{1-\gamma}{\gamma}}
\right\}.  
\end{equation}
\end{proposition}The result of Proposition~\ref{prop:harmonic-upper} immediately indicates an
$O_\gamma(r^{-1/2})$ upper bound for every fixed $\gamma\in(0,1)$. However, it is worth noting that, as either $\gamma\downarrow0$ or $\gamma\uparrow1$, the term $\max\left\{
8\sqrt{\frac{\gamma}{1-\gamma}},
2\sqrt{\frac{1-\gamma}{\gamma}}
\right\}$ will diverge such that the
resulting upper bound becomes \emph{uninformative}.

Recall that standard forward greedy achieves the tight
$1/2$-approximation ratio for monotone submodular maximization under a
matroid constraint~\cite{fisher1978analysis}. Moreover, whenever
$0<\gamma_1\le\gamma_2\le1$, every $\gamma_2$-weakly submodular function is
also $\gamma_1$-weakly submodular. It follows that
\begin{equation}\label{add_submodular_endpoint}
    R_{\gr}(\gamma,r)\le R_{\gr}(1,r)=\frac12
\qquad\text{for every }\gamma\in(0,1].
\end{equation}
This universal $1/2$ upper bound can help us to control the behavior near $\gamma=1$. To obtain a sharper bound when $\gamma$ is small, we next introduce a more finely calibrated core.

\subsection{Second core: a general power mean}
As in the preceding subsection, we consider the same partition matroid
constraint with a general power-mean core.
For $p\ge1$, we define

\begin{equation}\label{eq:power-mean}
\Psi_p(x,s)\triangleq \frac{xs}{(x^p+s^p)^{1/p}},
\qquad x,s\ge0.
\end{equation}
At the origin, we set $\Psi_p(0,0):=0$.
We first explore its continuous properties and then pass them to the induced set
function.  The proofs are deferred to Appendix~\ref{app:upper-proofs}.

\begin{lemma}[Continuous properties of the power mean]
\label{lem:power-continuous}
The function $\Psi_p$ is nonnegative, coordinate-wise nondecreasing,
concave, and bounded above by $\min\{x,s\}$.
For $x,s>0$,
\begin{align*}
\partial_x\Psi_p(x,s)=s^{p+1}(x^p+s^p)^{-1-1/p},\qquad
\partial_s\Psi_p(x,s)=x^{p+1}(x^p+s^p)^{-1-1/p}.
\end{align*}
Moreover, the right derivative $\partial_s^+\Psi_p(x,0)$ equals one for
every $x>0$.
\end{lemma}We next investigate a set objective induced by the aforementioned core of the power mean. Specifically, for
$a\ge1$, we define
\[
\widehat\Psi_{a,p}(S)
=\Psi_p\bigl(a|S\cap G|,|S\cap E|\bigr),
\qquad S\subseteq U,
\] where $G$ denotes the set of \emph{good} elements and $E$ an \emph{auxiliary} set. Subsequently, we establish the weak submodularity of $\widehat\Psi_{a,p}$ when the starting set
contains at least one good element, that is,

\begin{lemma}
\label{lem:power-ratio}
For every inclusion $A\subseteq C\subseteq U$ with $|A\cap G|\ge1$,
\begin{equation}\label{eq:power-discrete-ratio}
\sum_{e\in C\setminus A}
\bigl(\widehat\Psi_{a,p}(A+e)-\widehat\Psi_{a,p}(A)\bigr)
\ge c_p\bigl(\widehat\Psi_{a,p}(C)-\widehat\Psi_{a,p}(A)\bigr),
\end{equation}
where $c_p\triangleq(1-2^{-p})/p$.
\end{lemma}Next, we consider the following special parameterization:
\begin{equation}\label{eq:tau-equation}
\gamma\triangleq
\frac{\tau}{(1+\tau)(1+\tau^p)^{1/p}},
\end{equation}
where $\tau\in(0,1]$ and let  $0<\gamma\le\min\{c_p,2^{-1-1/p}\}$. Note that this requirement is possible because the right-hand side of
Eq.\eqref{eq:tau-equation}, viewed as a function of $\tau$, increases
continuously from zero to $2^{-1-1/p}$.  Thus, such a suitable $\tau$ exists. Then, for an integer $k\ge1$, we set $a=k/\tau\ge1$ and define the
corresponding count core as follows:
\begin{equation}\label{eq:capped-core}
H_{p,\tau,k}(q,s)=
\begin{cases}
\min\{s,k\},&q=0,\\[1mm]
\gamma^{-1}\Psi_p(aq,s),&q\ge1.
\end{cases}
\end{equation}
From Eq.~\ref{eq:capped-core}, we know this core $H_{p,\tau,k}$ controls the behavior that crosses from zero good elements to one or more
good elements via the interaction between $\Psi_p$ and $\min\{s,k\}$. The next lemma establishes the weak submodularity of this core $H_{p,\tau,k}$, namely,

\begin{lemma}
\label{lem:capped-valid}
Under Eq.\eqref{eq:tau-equation} and $0<\gamma\le\min\{c_p,2^{-1-1/p}\}$, we define
\[
\widehat H_{p,\tau,k}(S)
=H_{p,\tau,k}\bigl(|S\cap G|,|S\cap E|\bigr),
\qquad S\subseteq U.
\]
Then $\widehat H_{p,\tau,k}$ is normalized, monotone, and globally
$\gamma$-weakly submodular.
\end{lemma}

\subsubsection{Sharpening the joint dependence with the power core}

As in the first-core construction, we use the partition matroid introduced
at the beginning of of Section~\ref{sec:upper-bounds} and assume $N\ge k$.
We instantiate Eq.\eqref{eq:common-template} with $w\triangleq 1+\eps$ and the core
$H_{p,\tau,k}$, namely,
\begin{equation}\label{eq:capped-objective}
f(S)=(1+\eps)\cdot|S\cap B|+
H_{p,\tau,k}(|S\cap G|,|S\cap E|),
\qquad\eps>0.
\end{equation}
Here, $\gamma=\frac{\tau}{(1+\tau)(1+\tau^p)^{1/p}}$ and $0<\gamma\le\min\{c_p,2^{-1-1/p}\}$. Lemma~\ref{lem:capped-valid} shows that the set function induced by the core
$H_{p,\tau,k}$ is normalized, monotone, and $\gamma$-weakly submodular. Adding the nonnegative modular term preserves these properties, so the
objective $f$ in \eqref{eq:capped-objective} is also monotone $\gamma$-weakly submodular. When $S$ contains neither good nor auxiliary elements,
equivalently $S\subseteq B$, the relevant marginal values are
\[
f(b_i\mid S)=1+\eps,\qquad
f(g_i\mid S)=0,\qquad
f(e_j\mid S)=1.
\]
Thus the greedy method strictly chooses all $k$ bad elements first and
thereby blocks all good elements.  After this, the first $k$ auxiliary
elements have marginal value one, while each remaining auxiliary element
has marginal value zero.  The greedy method selects all $N$ auxiliary
elements and returns $S_r\triangleq B\cup E$.  Since $N\ge k$, substituting
$b(B\cup E)=k$, $q(B\cup E)=0$, and $s(B\cup E)=N$ gives
\begin{equation}\label{eq:capped-greedy-value}
f(S_r)=f(B\cup E)
=(2+\eps)k.
\end{equation}

Moreover, the base $G\cup E$ is feasible.  Since
$a=k/\tau$, direct substitution gives
\begin{equation}\label{eq:capped-good-value}
f(G\cup E)
=\frac{k^2N}
{\gamma\left(k^{2p}+\tau^pN^p\right)^{1/p}}.
\end{equation}
Comparing these two bases and letting $\eps\downarrow0$ yields
\begin{equation}\label{eq:capped-ratio}
R_{\gr}(\gamma,k+N)\le \frac{f(S_r)}{f(G\cup E)}
\le2\gamma
\left[
\left(\frac{k}{N}\right)^p+
\left(\frac{\tau}{k}\right)^p
\right]^{1/p}.
\end{equation}
As a result, from $r\triangleq N+k$ and $N\ge k$, we have 
\begin{proposition}
\label{prop:power-upper}For every $p\ge1$ and
$0<\gamma\le\min\{c_p,2^{-1-1/p}\}$, let $\tau\in(0,1]$ satisfy
Eq.\eqref{eq:tau-equation}.  Then, for every $r\ge2$,
\begin{equation}\label{eq:power-upper}
R_{\gr}(\gamma,r)
\le
\min_{1\le k\le\lfloor r/2\rfloor}
2\gamma
\left[
\left(\frac{k}{r-k}\right)^p+
\left(\frac{\tau}{k}\right)^p
\right]^{1/p}.
\end{equation}
\end{proposition}The main advantage of Proposition~\ref{prop:power-upper} mainly lies in how the weak submodularity ratio enters the resulting bound. Note that In
Eq~\ref{eq:power-upper}, $\gamma$ appears as a multiplicative
factor and does not occur in the denominator. Consequently, this feature allows us
to capture the correct behavior of $R_{\gr}(\gamma,r)$ in the small-$\gamma$
regime.

\subsection{A uniform upper bound}
The preceding subsections develop two complementary upper-bound
constructions. More specifically, Proposition~\ref{prop:power-upper} uses the
power core to recover the dependence on $\gamma$ that is missing as
$\gamma\downarrow0$. Once $\gamma$ is bounded away from zero, the simpler
harmonic construction of Proposition~\ref{prop:harmonic-upper} yields the
required $O_\gamma(r^{-1/2})$ dependence on the matroid rank. Finally, by
combining these two constructions with the classical tight
$1/2$-approximation result of submodular maximization,  we get the following theorem
\begin{theorem}\label{thm:uniform-upper}
For every $0<\gamma\le1$ and $r\ge2$, we can prove that
\begin{equation}\label{eq:uniform-upper}
R_{\gr}(\gamma,r)
\le\min\left\{\frac{1}{2},20\max\left\{\frac{\gamma}{r},
\frac{\gamma^2}{B_\gamma+
\sqrt{\gamma(1-\gamma)(r-1)}}\right\}\right\},
\end{equation} where $B_\gamma\triangleq2-\gamma+\gamma^2$.
\end{theorem}

\section{Lower bounds for
\texorpdfstring{$R_{\gr}(\gamma,r)$}{Rgr(gamma,r)}}
\label{sec:lower-bounds}
In this section,  we suppose that $O\in\cI$ is an optimal base, namely, 
$f(O)\triangleq\max_{T\in\cI}f(T)$, and set $\OPT\triangleq f(O)$.  For the greedy sets
$S_0,S_1,\ldots,S_r$, we define 
\[
F_t\triangleq f(S_t),\qquad d_t\triangleq F_t-F_{t-1},
\qquad
W_t\triangleq \sum_{i=1}^t i d_i,
\qquad
K_t\triangleq\sum_{j=0}^{t-1}F_j.
\]
We now turn to lower bounds for $R_{\gr}(\gamma,r)$. We first establish a
universal approximation guarantee for standard forward greedy that holds for
every $\gamma\in(0,1]$ and every matroid rank $r$. We then develop a refined
analysis for the regime in which $\gamma$ is close to one. Note that all proofs are collected in Appendix~\ref{app:lower-proofs}.

\subsection{A universal rank-dependent guarantee}
In this subsection, we begin with several lemmas that will be frequently used in the lower-bound analysis of $R_{\gr}(\gamma,r)$ and then state our obtained lower bound. More specifically, we first investigate the relationship between $F_{t}$, $W_t$ and $K_t$.
\begin{lemma}[Basic identities]\label{lem:basic-identities}
For every $1\le t\le r$, we have 
\begin{enumerate}[label=\textup{(\roman*)},leftmargin=2.2em]
\item $d_t\ge0$;
\item $W_t=tF_t-K_t$ and $W_t\ge F_t$.
\end{enumerate}
\end{lemma}

The next lemma orders the elements of $O$ so that every greedy solution $S_t$ can be
completed to a base using the remaining elements of $O$.

\begin{lemma}\label{lem:ordering}
The optimal base can be ordered as $O=\{o_1,\ldots,o_r\}$ so that, for any $t=0,\dots,r-1$,
\[
S_t\cup\{o_{t+1},\ldots,o_r\}
\]
is a base.
In particular, $o_i$ is feasible at $S_{i-1}$, and $o_j$ is feasible at
$S_t$ whenever $j>t$.
\end{lemma}With this result of Lemma~\ref{lem:ordering} and $\gamma$-weak submodularity, we also can show that
\begin{lemma}\label{lem:delayed}
For $1\le i\le t\le r$,
\begin{equation}\label{eq:delayed}
\gamma f(o_i\mid S_t)
\le (t-i+2)d_i-\gamma(F_t-F_{i-1}).
\end{equation}
\end{lemma}By summing Eq.~\eqref{eq:delayed}, we also can have
\begin{lemma}\label{lem:prefix}
Let $B_\gamma\triangleq 2-\gamma+\gamma^2$.  For every $0\le t\le r$,
\begin{equation}\label{eq:prefix-main}
\gamma^2\OPT
\le B_\gamma F_t+K_t+\gamma(r-t)d_{t+1}.
\end{equation}
Particularly when $t=r$, we also can show 
\begin{equation}\label{eq:prefix-terminal}
\gamma^2\OPT\le B_\gamma F_r+K_r.
\end{equation}
\end{lemma}To use the relations in Eq.~\eqref{eq:prefix-main} and Eq.~\eqref{eq:prefix-terminal}, we first establish the following recursive inequality:
\begin{lemma}[Recurrence aggregation]\label{lem:recurrence}
Let $r\ge1$ and $A,B,C\ge0$, let
$0=F_0\le F_1\le\cdots\le F_r=M$, and set
$K_t=\sum_{j<t}F_j$.  Suppose that
\[
C\le A(F_{t+1}-F_t)+BF_t+K_t\quad(t<r),
\qquad C\le BM+K_r.
\]
Then
\begin{equation}\label{eq:recurrence-conclusion}
C\le\left(B+\sqrt{A(1-1/r)}\right)M.
\end{equation}
\end{lemma}
According to this Lemma~\ref{lem:recurrence} and the results of Lemma~\ref{lem:delayed}, we can prove \begin{theorem}\label{thm:lower}
For every $0<\gamma\le1$ and $r\ge2$, if $B_\gamma\triangleq 2-\gamma+\gamma^2$, 
\begin{equation}\label{eq:main-lower}
R_{\gr}(\gamma,r)
\ge
\max\left\{\frac{\gamma}{r},
\frac{\gamma^2}{B_\gamma+\sqrt{\gamma(r-1)}}
\right\}.
\end{equation}
\end{theorem}It is worth noting that, as $\gamma \rightarrow 1$, the right-hand side of
Eq.\eqref{eq:main-lower} will converge to
\[
\max\left\{\frac{1}{r},\frac{1}{2+\sqrt{r-1}}\right\}.
\]
For $r\ge 3$, this bound remains rank-dependent and is \emph{strictly} weaker than the classical $1/2$-approximation guarantee for monotone submodular maximization~\cite{fisher1978analysis}. A natural question is \emph{whether this gap can be closed when $\gamma=1$?}. We next show that this is indeed possible through an endpoint correction, that is to say, \begin{theorem}\label{thm:endpoint}
For every $0<\gamma\le1$ and $r\ge1$, we can show that
\begin{equation}\label{eq:endpoint-lower}
R_{\gr}(\gamma,r)
\ge
\frac{\gamma^{r+1}}{2-\gamma+\gamma^{r+1}}.
\end{equation}
\end{theorem} Note that $\lim_{\gamma\uparrow1}\frac{\gamma^{r+1}}{2-\gamma+\gamma^{r+1}}=\frac{1}{2}$. Then, combining both Theorem~\ref{thm:lower} and Theorem~\ref{thm:endpoint}, we can obtain that, for every $0<\gamma\le 1$ and $r\ge 2$,
\begin{equation}\label{add_1}
  R_{\gr}(\gamma,r)
\ge
\max\left\{
\frac{\gamma}{r},
\frac{\gamma^2}{B_\gamma+\sqrt{\gamma(r-1)}},
\frac{\gamma^{r+1}}{2-\gamma+\gamma^{r+1}}
\right\},  
\end{equation}
where $B_\gamma\triangleq2-\gamma+\gamma^2$. We close this subsection by comparing our obtained lower bound in Eq.~\eqref{add_1} with the two lower guarantees from the analysis of \cite{gatmiry2018nonsubmodular}.\begin{lemma}[Comparison with previous lower bounds]
\label{lem:lower-comparison}
For every $0<\gamma\le1$ and $r\ge2$,
\begin{equation}\label{eq:lower-comparison}
R_{\gr}(\gamma,r)\ge\frac{\gamma^2}{B_\gamma+\sqrt{\gamma(r-1)}}
>
\frac{\gamma^2}{2+2\sqrt{\gamma r}+\gamma^2}
>
\frac{0.4\gamma^2}{1+\sqrt{\gamma r}}.
\end{equation}
Thus, for $r\ge2$, our obtained lower bound in Eq.\eqref{add_1} strictly improves the
earlier result in \cite{gatmiry2018nonsubmodular}.
\end{lemma}
\subsection{A Refined lower-bound analysis for \texorpdfstring{$\gamma\in[1/2,1)$}{}}
Recall that Theorem~\ref{thm:lower} gives
\[R_{\gr}(\gamma,r)
\ge
\max\left\{
\frac{\gamma}{r},
\frac{\gamma^2}
{B_\gamma+\sqrt{\gamma(r-1)}}
\right\},\] whereas the upper bound in Theorem~\ref{thm:uniform-upper} is governed by
\[
\Phi(\gamma,r)
\triangleq
\max\left\{
\frac{\gamma}{r},
\frac{\gamma^2}
{B_\gamma+\sqrt{\gamma(1-\gamma)(r-1)}}
\right\}.
\]
When $0<\gamma\le1/2$, we can show these two expressions match up to an
absolute constant factor. More specifically,  if $1-\gamma\ge1/2$, we have $\sqrt{\gamma(r-1)}
\le
\sqrt{2\gamma(1-\gamma)(r-1)}$ such that 
\[
R_{\gr}(\gamma,r)
\ge
\max\left\{
\frac{\gamma}{r},
\frac{\gamma^2}
{\sqrt2\left(B_\gamma+\sqrt{\gamma(1-\gamma)(r-1)}\right)}
\right\}.
\]
As a result, over this range $0<\gamma\le1/2$, our obtained  universal lower bound in Theorem~\ref{thm:lower} captures the same 
asymptotic dependence on $\gamma$ and $r$ as the upper bound of Theorem~\ref{thm:uniform-upper}.

However, this relationship breaks down as
$\gamma\uparrow1$. The
fundamental reason is that the universal lower bound contains
$\sqrt{\gamma(r-1)}$ in its denominator, whereas the upper profile contains
the smaller term $\sqrt{\gamma(1-\gamma)(r-1)}$. To overcome this challenge, we next develop a refined lower-bound analysis for $\gamma\in[1/2,1)$, which can explicitly preserves the
missing factor $1-\gamma$.

Before going into the details, like Lemma~\ref{lem:prefix} we first present a sharper upper bound on $\OPT$ for the regime where $\gamma$ is close to one. More specifically, we have 
\begin{lemma}
\label{lem:near-prefix}
For $1/2\le\gamma<1$ and every $0\le t<r$, if $D_\gamma\triangleq60-15\gamma+\gamma^2$, we have
\begin{equation}\label{eq:near-prefix-main}
\gamma^2\OPT
\le D_\gamma F_t+15(1-\gamma)K_t+\gamma(r-t)d_{t+1}.
\end{equation}
In particular,
\begin{equation}\label{eq:near-prefix-terminal}
\gamma^2\OPT\le D_\gamma F_r+15(1-\gamma)K_r.
\end{equation}
\end{lemma}
Recall that
Lemma~\ref{lem:prefix} gives the result $\gamma^2\OPT
\le
B_\gamma F_t+K_t+\gamma(r-t)d_{t+1}$. In contrast, this new Lemma~\ref{lem:near-prefix} replaces the coefficient of the accumulated
history $K_t$ by $15(1-\gamma)$. Then, applying the same recurrence argument to Eq.\eqref{eq:near-prefix-main} and  Eq.\eqref{eq:near-prefix-terminal} can produce a better dependence of order $\sqrt{\gamma(1-\gamma)(r-1)}$, that is to say,

\begin{theorem}\label{thm:near-lower}
For every $1/2\le\gamma<1$ and $r\ge2$, if $D_\gamma\triangleq60-15\gamma+\gamma^2$, we can show that 
\begin{equation}\label{eq:near-lower}
R_{\gr}(\gamma,r)
\ge
\max\left\{\frac{\gamma}{r},
\frac{\gamma^2}
{D_\gamma+\sqrt{15\gamma(1-\gamma)(r-1)}}\right\}.
\end{equation}
\end{theorem}
Merging the result of Theorem~\ref{thm:near-lower} into the aforementioned Eq.\eqref{add_1}, we can finally obtain a uniform lower bound for $R_{\gr}(\gamma,r)$, namely,
\begin{theorem}
\label{thm:uniform-lower}
For every $0<\gamma\le1$ and every integer $r\ge2$, we can prove that
\begin{equation}\label{eq:uniform-lower}
R_{\gr}(\gamma,r)\ge \max\left\{
\frac{\gamma}{r},
\frac{\gamma^2}
{B_\gamma+\sqrt{\gamma(r-1)}},
\frac{\gamma^{r+1}}
{2-\gamma+\gamma^{r+1}},
\frac{\gamma^2}
{14\left(B_\gamma+\sqrt{\gamma(1-\gamma)(r-1)}\right)}
\right\},
\end{equation} where $B_\gamma\triangleq2-\gamma+\gamma^2$. 
\end{theorem}

\section{Conclusion}

In this paper, we provide an asymptotically tight characterization of the
worst-case performance of standard forward greedy for matroid-constrained
weakly submodular maximization. Specifically, for every submodularity ratio
$\gamma\in(0,1]$ and matroid rank $r\ge2$, we establish
\[
R_{\mathrm{gr}}(\gamma,r)
=
\Theta\!\left(
\max\left\{
\frac{\gamma}{r},
\frac{\gamma^2}
{1+\sqrt{\gamma(1-\gamma)r}}
\right\}
\right).
\]
This characterization reveals how the performance of standard greedy varies
across different $\gamma$ and $r$. In particular, when $\gamma\to0$ and $\gamma r=O(1)$, the first term in the
maximum will dominate, and hence $R_{\mathrm{gr}}(\gamma,r)
=
\Theta\!\left(\gamma/r\right)$. By contrast, if $\gamma\to0$ while $\gamma r\to\infty$, then the second term
becomes dominant and yields $R_{\mathrm{gr}}(\gamma,r)
=
\Theta\!\left(\gamma^{3/2}/\sqrt{r}
\right)$. Moreover, for every fixed $\gamma\in(0,1)$, our characterization gives $R_{\mathrm{gr}}(\gamma,r)
=
\Theta_{\gamma}\!\left(r^{-1/2}\right)$, which indicates that even when $\gamma$ is fixed, the worst-case approximation guarantee of standard greedy will deteriorate at the rate $r^{-1/2}$ as the matroid rank grows. Overall, our result provides a complete description of how weak
submodularity and the matroid rank jointly influence the worst-case
performance of standard forward greedy.

\bibliographystyle{plainnat}
\bibliography{references}
\appendix

\section*{Contents of the Appendices}
\phantomsection
\label{app:contents}
\appendixcontents

\section{Proofs for the upper bounds}
\label{app:upper-proofs}
\addappendixsection{Proofs for the upper bounds}
\subsection{Proof of Lemma~\ref{lem:H-continuous}}
\label{app:proof-H-continuous}
\addappendixsubsection{Proof of Lemma~\ref{lem:H-continuous}}
\begin{proof}
For the harmonic mean $\phi(x,s)=2xs/(x+s)$, we can show its inner matrix is positive semidefinite, namely,
\[
\nabla^2\phi(x,s)
=-\frac{4}{(x+s)^3}
\begin{pmatrix}s^2&-xs\\-xs&x^2\end{pmatrix}
\preceq0.
\]
As a result, $\phi(x,s)$ is concave. Generally speaking, concavity can be
preserved under affine transformation and positive scaling, so $H(q,s)=\frac12\phi(a(q+\lambda),s)$ is also concave.  Direct differentiation gives
\[
\partial_q H(q,s)=\frac{as^2}{D^2},
\qquad
\partial_s H(q,s)=\frac{x^2}{D^2},
\] where $x\triangleq a(q+\lambda)$ and $D\triangleq x+s$. Importantly, from the aforementioned two identities, we can easily know that the partial derivatives are nonnegative, so $H(q,s)$ is coordinate-wise nondecreasing. Finally, based on the definition of $H(q,s)$, we can also have that
\begin{align*}
H(q+1,s)-H(q,s)
&=\frac{as^2}{D(D+a)}
=\frac{D}{D+a}\partial_q H(q,s),\\
H(q,s+1)-H(q,s)
&=\frac{x^2}{D(D+1)}
=\frac{D}{D+1}\partial_s H(q,s),
\end{align*}
which proves the claimed identities in Lemma~\ref{lem:H-continuous}.
\end{proof}

\subsection{Proof of Lemma~\ref{lem:H-ratio}}
\label{app:proof-H-ratio}
\addappendixsubsection{Proof of Lemma~\ref{lem:H-ratio}}
\begin{proof}
Normalization follows from
$\widehat H(\emptyset)=H(0,0)=\frac{1}{2}\phi(a\lambda,0)=0$. Moreover, from  Lemma~\ref{lem:H-continuous}, we know that $H$ is coordinate-wise nondecreasing, so $\widehat H(S)\triangleq H\bigl(|S\cap G|,|S\cap E|\bigr)$  is monotone. Fix
$A\subseteq C\subseteq U$ and put
\[
q=|A\cap G|,\quad s=|A\cap E|,
\quad u=|(C\setminus A)\cap G|,
\quad v=|(C\setminus A)\cap E|.
\]
Elements of $B$ have zero marginal under $\widehat H$, and hence
\begin{align}
\sum_{e\in C\setminus A}\widehat H(e\mid A)
&=u\bigl(H(q+1,s)-H(q,s)\bigr)
+v\bigl(H(q,s+1)-H(q,s)\bigr).
\label{eq:H-singleton-sum}
\end{align}
With $D\triangleq a(q+\lambda)+s$ and $a\ge1$, we have $D\ge a\lambda$.
Therefore
\begin{align*}
\frac{D}{D+a}&\ge\frac{\lambda}{1+\lambda},&
\frac{D}{D+1}&\ge\frac{\lambda}{1+\lambda}.
\end{align*}
Together with Lemma~\ref{lem:H-continuous} and $\gamma_{\lambda}\triangleq\frac{\lambda}{1+\lambda}$, these inequalities imply
\begin{align*}
H(q+1,s)-H(q,s)
&\ge\gamma_\lambda\cdot\partial_q H(q,s),\\
H(q,s+1)-H(q,s)
&\ge\gamma_\lambda\cdot\partial_s H(q,s).
\end{align*} Using Eq.\eqref{eq:H-singleton-sum} and then the concavity of $H$, we obtain
\begin{align*}
\sum_{e\in C\setminus A}\widehat H(e\mid A)
&\ge\gamma_\lambda
\bigl(u\partial_q H(q,s)+v\partial_s H(q,s)\bigr)\\
&\ge\gamma_\lambda
\bigl(H(q+u,s+v)-H(q,s)\bigr)\\
&=\gamma_\lambda\bigl(\widehat H(C)-\widehat H(A)\bigr).
\end{align*}
So $\widehat H(S)\triangleq H\bigl(|S\cap G|,|S\cap E|\bigr)$  is $\gamma_\lambda$-weakly submodular.
\end{proof}

\subsection{Proof of Proposition~\ref{prop:harmonic-upper}}
\label{app:proof-harmonic-upper}
\addappendixsubsection{Proof of Proposition~\ref{prop:harmonic-upper}}
\begin{proof}
Note that in Section~\ref{sec:how_the_harmonic}, we assume $N=r-k\ge k$ such that
$1\le k\le\lfloor r/2\rfloor$.  Combining this observation with
Eq.\eqref{eq:offset-ratio}, we thus have
\begin{equation*}
    R_{\gr}(\gamma,r)
\le
\min_{1\le k\le\lfloor r/2\rfloor}
\frac{2k+\lambda}{k+\lambda}
\left(\frac{k}{r-k}+\frac{\lambda}{k+\lambda}\right),
\end{equation*} where $\lambda\triangleq\frac{\gamma}{1-\gamma}$. Next, we show that
\[
R_{\gr}(\gamma,r)
\le\frac1{\sqrt r}
\max\left\{8\sqrt\lambda,\frac{2}{\sqrt\lambda}\right\}.
\]
First, we suppose that
$\lambda r\ge4$ and $r\ge4\lambda$, and set
\[
k=\left\lfloor\sqrt{\lambda r}\right\rfloor,
\qquad N=r-k.
\]
Then
\[
\frac12\sqrt{\lambda r}\le k\le\sqrt{\lambda r}\le\frac r2,
\qquad N=r-k\ge\frac r2,
\]
so Eq.\eqref{eq:offset-ratio} gives
\[
R_{\gr}(\gamma,r)\le 2\left(\frac{k}{N}+\frac{\lambda}{k}\right)
\le2\left(2\sqrt{\frac{\lambda}{r}}
          +2\sqrt{\frac{\lambda}{r}}\right)
=\frac{8\sqrt\lambda}{\sqrt r}
\le\frac1{\sqrt r}
\max\left\{8\sqrt\lambda,\frac{2}{\sqrt\lambda}\right\}.
\]
If $\lambda r<4$, then
\[
\frac1{\sqrt r}
\max\left\{8\sqrt\lambda,\frac{2}{\sqrt\lambda}\right\}
\ge \frac{2}{\sqrt{\lambda r}}
>1
\ge R_{\gr}(\gamma,r).
\]
If $r<4\lambda$, then
\[
\frac1{\sqrt r}
\max\left\{8\sqrt\lambda,\frac{2}{\sqrt\lambda}\right\}
\ge \frac{2\sqrt\lambda}{\sqrt r}
>1
\ge R_{\gr}(\gamma,r).
\]
Thus, in all cases, we obtain
\[
R_{\gr}(\gamma,r)
\le\frac1{\sqrt r}
\max\left\{8\sqrt\lambda,\frac{2}{\sqrt\lambda}\right\}.
\]

\end{proof}

\subsection{Proof of Lemma~\ref{lem:power-continuous}}
\label{app:proof-power-continuous}
\addappendixsubsection{Proof of Lemma~\ref{lem:power-continuous}}
\begin{proof}
Set $Q\triangleq x^p+s^p$.  For $x,s>0$, direct differentiation gives
\begin{align*}
\partial_x\Psi_p(x,s)&=s^{p+1}Q^{-1-1/p},&
\partial_s\Psi_p(x,s)&=x^{p+1}Q^{-1-1/p}.
\end{align*}
Both partial derivatives are nonnegative, so $\Psi_p$ is coordinate-wise nondecreasing. Moreover, from the definition of $\Psi_p$, we can know 
the Hessian is
\[
\nabla^2\Psi_p(x,s)
=-(p+1)x^{p-1}s^{p-1}Q^{-2-1/p}
\begin{pmatrix}
s^2&-xs\\
-xs&x^2
\end{pmatrix}.
\]
Its diagonal entries are nonpositive and its determinant is zero; hence the
Hessian is negative semidefinite, and $\Psi_p$ is concave on the positive
orthant. As for the concavity over boundary points, we can apply the concave inequality to positive sequences and let these positive sequences approximate the boundary points. In addition, for any $x>0$,
\[
\partial_s^+\Psi_p(x,0)
=\lim_{h\downarrow0}
\frac{\Psi_p(x,h)-\Psi_p(x,0)}{h}
=\lim_{h\downarrow0}\frac{x}{(x^p+h^p)^{1/p}}=1.
\] Finally, the inequality $(x^p+s^p)^{1/p}\ge\max\{x,s\}$ gives
$\Psi_p(x,s)\le\min\{x,s\}$.
\end{proof}

\subsection{Proof of Lemma~\ref{lem:power-ratio}}
\label{app:proof-power-ratio}
\addappendixsubsection{Proof of Lemma~\ref{lem:power-ratio}}
\begin{proof}
Fix $A\subseteq C\subseteq U$, and set
\[
q=|A\cap G|,\qquad s=|A\cap E|,
\qquad u=|(C\setminus A)\cap G|,
\qquad v=|(C\setminus A)\cap E|.
\]
Since $q\ge1$, the quantity $x=aq$ satisfies $x\ge a$.  Elements of $B$
have zero marginal under $\widehat\Psi_{a,p}$, and therefore
\begin{align}
&\sum_{e\in C\setminus A}
\bigl(\widehat\Psi_{a,p}(A+e)-\widehat\Psi_{a,p}(A)\bigr)\notag\\
&\quad =u\bigl(\Psi_p(x+a,s)-\Psi_p(x,s)\bigr)
+v\bigl(\Psi_p(x,s+1)-\Psi_p(x,s)\bigr).
\label{eq:power-singleton-sum}
\end{align}

We first consider the marginal of a good element.  If $s=0$, then both
$\Psi_p(x+a,s)-\Psi_p(x,s)$ and $\partial_x\Psi_p(x,s)$ are zero.  If
$s>0$, the fundamental theorem of calculus gives
\begin{align*}
\Psi_p(x+a,s)-\Psi_p(x,s)
&=a\int_0^1\partial_x\Psi_p(x+at,s)\,dt.
\end{align*}
For every $t\in[0,1]$, the inequality $x+at\le x(1+t)$ implies
\begin{align*}
\frac{\partial_x\Psi_p(x+at,s)}{\partial_x\Psi_p(x,s)}=\left(\frac{x^p+s^p}{(x+at)^p+s^p}\right)^{1+1/p}\ge (1+t)^{-p-1}.
\end{align*}
Consequently, for any $s\ge0$, we have 
\begin{align}
\Psi_p(x+a,s)-\Psi_p(x,s)\ge a\partial_x\Psi_p(x,s)
\int_0^1(1+t)^{-p-1}\,dt\notag=c_p a\partial_x\Psi_p(x,s),
\label{eq:power-forward-x}
\end{align} where $c_p\triangleq(1-2^{-p})/p$. We next consider the marginal of an auxiliary element.  When $s\ge1$,
the inequality $s+t\le s(1+t)$ similarly gives
\begin{align*}
\frac{\partial_s\Psi_p(x,s+t)}{\partial_s\Psi_p(x,s)}=\left(\frac{x^p+s^p}{x^p+(s+t)^p}\right)^{1+1/p}\ge(1+t)^{-p-1}.
\end{align*}
Thus, for $s\ge1$, we have 
\begin{equation}
\Psi_p(x,s+1)-\Psi_p(x,s)
\ge c_p\partial_s\Psi_p(x,s).
\label{eq:power-forward-s}
\end{equation}When $s=0$, since $\Psi_p(x,0)=0$, we have
\begin{align*}
\Psi_p(x,1)-\Psi_p(x,0)=\Psi_p(x,1)=\frac{x}{(x^p+1)^{1/p}}\ge2^{-1/p}\ge c_p
=c_p\partial_s^+\Psi_p(x,0),
\end{align*}
where the first inequality comes from $x\ge a\ge1$ and
the second inequality follows from $p\ge1$and $c_p\le1/2\le2^{-1/p}$. It follows from Eq.\eqref{eq:power-singleton-sum}-Eq.\eqref{eq:power-forward-s}
that
\begin{align}
\sum_{e\in C\setminus A}
\bigl(\widehat\Psi_{a,p}(A+e)-\widehat\Psi_{a,p}(A)\bigr)\ge c_p\bigl(au\partial_x\Psi_p(x,s)
+v\partial_s\Psi_p(x,s)\bigr),
\label{eq:power-gradient-lower}
\end{align}where $\partial_s$ is interpreted as the right derivative when $s=0$. By the concavity of $\Psi_p$, for any $s\ge 1$
\begin{align*}
\Psi_p(x+au,s+v)-\Psi_p(x,s)
\le au\partial_x\Psi_p(x,s)+v\partial_s\Psi_p(x,s).
\end{align*} As for $s=0$, a similar inequality also follows directly from
$\Psi_p(x+au,v)\le v$, together with
$\partial_x\Psi_p(x,0)=0$ and $\partial_s^+\Psi_p(x,0)=1$. More specifically, we have 
\begin{align*}
\Psi_p(x+au,v)-\Psi_p(x,0)=\Psi_p(x+au,v)\le v=au\partial_x\Psi_p(x,0)+v\partial_s^+\Psi_p(x,0).
\end{align*}Combining this inequality with Eq.\eqref{eq:power-gradient-lower} gives
Eq.\eqref{eq:power-discrete-ratio}.
\end{proof}

\subsection{Proof of Lemma~\ref{lem:capped-valid}}
\label{app:proof-capped-valid}
\addappendixsubsection{Proof of Lemma~\ref{lem:capped-valid}}
\begin{proof}
Since $H_{p,\tau,k}(0,0)=0$, the induced set function is normalized.  To
prove the monotonicity and weak submodularity, we fix $A\subseteq C\subseteq U$, and
set
\[
q=|A\cap G|,\qquad s=|A\cap E|,
\qquad u=|(C\setminus A)\cap G|,
\qquad v=|(C\setminus A)\cap E|.
\]
Elements of $B$ have zero marginal under $\widehat H_{p,\tau,k}$.  Hence, by
the definition of $\widehat H_{p,\tau,k}$, the sum of the singleton marginals
at $A$ is
\begin{align*}
\sum_{e\in C\setminus A}\widehat H_{p,\tau,k}(e\mid A)=u\bigl(H_{p,\tau,k}(q+1,s)-H_{p,\tau,k}(q,s)\bigr)+v\bigl(H_{p,\tau,k}(q,s+1)-H_{p,\tau,k}(q,s)\bigr),
\end{align*}
and the joint gain is $H_{p,\tau,k}(q+u,s+v)-H_{p,\tau,k}(q,s)$.

We verify the weak submodularity in three cases.  First, we suppose
that $q\ge1$. Then, both the starting and ending counts then lie in the second
branch of Eq.\eqref{eq:capped-core}. Lemma~\ref{lem:power-ratio} and the
inequality $\gamma\le c_p$ immediately give
\begin{align*}
\sum_{e\in C\setminus A}\widehat H_{p,\tau,k}(e\mid A)\ge c_p\bigl(H_{p,\tau,k}(q+u,s+v)
                 -H_{p,\tau,k}(q,s)\bigr)\ge \gamma\bigl(H_{p,\tau,k}(q+u,s+v)
                    -H_{p,\tau,k}(q,s)\bigr).
\end{align*}

Second, suppose that $q=u=0$.  If $s\ge k$, then we have
\begin{equation*}
    \sum_{e\in C\setminus A}\widehat H_{p,\tau,k}(e\mid A)=0=H_{p,\tau,k}(q+u,s+v)-H_{p,\tau,k}(q,s).
\end{equation*}If $s<k$, every new auxiliary element has singleton
marginal one, and hence
\begin{align*}
&\sum_{e\in C\setminus A}\widehat H_{p,\tau,k}(e\mid A)=v\bigl(H_{p,\tau,k}(q,s+1)-H_{p,\tau,k}(q,s)\bigr)=v,\\
&H_{p,\tau,k}(0,s+v)-H_{p,\tau,k}(0,s)
=\min\{s+v,k\}-s\le v.
\end{align*}
Since $\gamma\le1$,  $\sum_{e\in C\setminus A}\widehat H_{p,\tau,k}(e\mid A)\ge\gamma(H_{p,\tau,k}(0,s+v)-H_{p,\tau,k}(0,s))$.

It remains to consider $q=0$ and $u\ge1$.  Define the singleton marginals
of a good element and an auxiliary element, respectively, by
\begin{align*}
\Delta_g(s)=H_{p,\tau,k}(1,s)-H_{p,\tau,k}(0,s),\quad \Delta_e(s)=H_{p,\tau,k}(0,s+1)-H_{p,\tau,k}(0,s).
\end{align*}
Then, we have  $\sum_{e\in C\setminus A}\widehat H_{p,\tau,k}(e\mid A)=u\Delta_g(s)+v\Delta_e(s)$.

Suppose first that $s<k$.  In this case $\Delta_e(s)=1$ and
$s/a=(s/k)\tau\le\tau$.  For $s>0$, from $\gamma\triangleq
\frac{\tau}{(1+\tau)(1+\tau^p)^{1/p}}$, we can show that
\begin{align*}
\gamma^{-1}\Psi_p(a,s)=\frac{\gamma^{-1}s}{(1+(s/a)^p)^{1/p}}\ge\frac{\gamma^{-1}s}{(1+\tau^p)^{1/p}}=\frac{1+\tau}{\tau}s\ge2s.
\end{align*}
Therefore, we can show that $\Delta_g(s)=\gamma^{-1}\Psi_p(a,s)-s\ge s$. The same inequality holds when $s=0$, because both sides are zero.  Using
$u\ge1$, we obtain
\begin{equation*}
\sum_{e\in C\setminus A}\widehat H_{p,\tau,k}(e\mid A)=u\Delta_g(s)+v\Delta_e(s)
\ge us+v\ge s+v.
\end{equation*}
On the other hand,
\[
\gamma\bigl(H_{p,\tau,k}(u,s+v)-H_{p,\tau,k}(0,s)\bigr)
=\Psi_p(au,s+v)-\gamma s
\le\Psi_p(au,s+v)\le s+v.
\]
Combining the last two displays proves the $\gamma$-weak submodularity when $s<k$.

Finally, we suppose that $s\ge k$.  Then $\Delta_e(s)=0$.  Since
$\Psi_p(a,s)$ is nondecreasing in $s$, so is $\Delta_g(s)=H_{p,\tau,k}(1,s)-k$.  Moreover,
$k=a\tau$ and Eq.\eqref{eq:tau-equation} give
\begin{align*}
\Delta_g(k)=\gamma^{-1}\Psi_p(a,a\tau)-a\tau=\frac{a\tau}{\gamma(1+\tau^p)^{1/p}}-a\tau=a(1+\tau)-a\tau=a.
\end{align*}
Consequently, $\sum_{e\in C\setminus A}\widehat H_{p,\tau,k}(e\mid A)=u\Delta_g(s)\ge u\Delta_g(k)= ua$.

The joint gain satisfies
\begin{align*}
\gamma\bigl(H_{p,\tau,k}(u,s+v)-H_{p,\tau,k}(0,s)\bigr)=\Psi_p(au,s+v)-\gamma k
\le\Psi_p(au,s+v)\le au.
\end{align*}
Thus the weak submodularity also holds when $s\ge k$.
\end{proof}

\subsection{Proof of Proposition~\ref{prop:power-upper}}
\label{app:proof-power-upper}
\addappendixsubsection{Proof of Proposition~\ref{prop:power-upper}}
\begin{proof}
For the power-mean upper bound, we assume $N=r-k\ge k$ such that
$1\le k\le\lfloor r/2\rfloor$.  Combining this observation with
Eq.\eqref{eq:capped-ratio}, we then can get the upper bound in
Eq.\eqref{eq:power-upper}.
\end{proof}

\subsection{Proof of Theorem~\ref{thm:uniform-upper}}
\label{app:proof-uniform-upper}
\addappendixsubsection{Proof of Theorem~\ref{thm:uniform-upper}}
\begin{proof}
In this section, we prove the uniform upper bound presented in Eq.\eqref{eq:uniform-upper}.  We divide the parameter range
into three parts.  First, we suppose that
$0<\gamma\le1/(2\sqrt2)$ and take $p=2$ in Proposition~\ref{prop:power-upper}. Moreover, from the definition of Eq.\eqref{eq:tau-equation}, namely $\gamma\triangleq\frac{\tau}{(1+\tau)\sqrt{1+\tau^2}}$, we can immediately get that $\frac{\tau}{2\sqrt{2}}\le\gamma\le\tau$. If $\tau r<4$, we choose $k=1$ in Eq.\eqref{eq:power-upper} such that
\begin{align*}
R_{\gr}(\gamma,r)\le2\gamma
\left[\frac1{(r-1)^2}+\tau^2\right]^{1/2}\le2\gamma
\left[\frac4{r^2}+\frac{16}{r^2}\right]^{1/2}
=4\sqrt5\,\frac{\gamma}{r}
<9\frac{\gamma}{r},
\end{align*} where the second inequality comes from $1/(r-1)\le2/r$($r\ge2$) and $\tau<4/r$.
Consequently, $R_{\gr}(\gamma,r)<9\frac{\gamma}{r}\le9\Phi(\gamma,r)$ where $\Phi(\gamma,r)\triangleq\max\left\{\frac{\gamma}{r},
\frac{\gamma^2}{B_\gamma+
\sqrt{\gamma(1-\gamma)(r-1)}}\right\}$. If $\tau r\ge4$, choose $k=\lfloor\sqrt{\tau r}\rfloor$.  Since
$\tau\le1$ and $\tau r\ge4$, we have $r\ge4$ and
\[
\frac12\sqrt{\tau r}\le k\le\sqrt{\tau r}
\le\sqrt r\le\frac r2,
\qquad r-k\ge\frac r2.
\]
Substitution into Eq.\eqref{eq:power-upper} gives
\begin{align}
R_{\gr}(\gamma,r)\le2\gamma
\left[
\left(\frac{k}{r-k}\right)^2+
\left(\frac{\tau}{k}\right)^2
\right]^{1/2}\le 2\gamma
\left[4\frac{\tau}{r}+4\frac{\tau}{r}\right]^{1/2}=4\sqrt2\,\gamma\sqrt{\frac{\tau}{r}}
\le2^{13/4}\frac{\gamma^{3/2}}{\sqrt r},
\label{eq:small-uniform-upper} 
\end{align} where the final inequality comes from $\tau\le2\sqrt{2}\gamma$. Let $y\triangleq \sqrt{\gamma r}$.  Since $B_\gamma\triangleq2-\gamma+\gamma^2\le2$ and
$\sqrt{\gamma(1-\gamma)(r-1)}\le y$, we thus have that
\[
\Phi(\gamma,r)\triangleq\max\left\{\frac{\gamma}{r},
\frac{\gamma^2}{B_\gamma+
\sqrt{\gamma(1-\gamma)(r-1)}}\right\}
\ge\gamma^2\max\{y^{-2},(2+y)^{-1}\}.
\]
If $y\le2$, then $\gamma^2/y\le2\gamma^2/y^2$; if $y\ge2$, then
$\gamma^2/y\le2\gamma^2/(2+y)$.  Hence, $\frac{\gamma^{3/2}}{\sqrt r}
=\frac{\gamma^2}{y}
\le2\Phi(\gamma,r)$. Together with Eq.\eqref{eq:small-uniform-upper}, this gives
$R_{\gr}(\gamma,r)\le2^{17/4}\Phi(\gamma,r)<20\Phi(\gamma,r)$.

Next, we suppose 
\[
\frac1{2\sqrt2}\le\gamma\le\frac12,
\qquad
\lambda\triangleq\frac{\gamma}{1-\gamma},
\qquad
x\triangleq\sqrt{\lambda r}.
\]
If $x<3$,  we have
\[
B_\gamma+
\sqrt{\gamma(1-\gamma)(r-1)}=B_\gamma+
(1-\gamma)\cdot \sqrt{\frac{\gamma}{1-\gamma}(r-1)}=B_\gamma+
(1-\gamma)x<B_\gamma+3(1-\gamma)\le 40\gamma^{2},
\] where the final inequality follows from that $B_\gamma+3(1-\gamma)=5-4\gamma+\gamma^2$ decreases over this
interval $\frac1{2\sqrt2}\le\gamma\le\frac12$, while $40\gamma^2$ increases, and at
$\gamma=1/(2\sqrt2)$, $B_\gamma+3(1-\gamma)=41/8-\sqrt2\le 5=40\gamma^2$.  As a result, according to Eq.\eqref{add_submodular_endpoint}, we have 
\[
R_{\gr}(\gamma,r)\le\frac12=20\frac{\gamma^{2}}{40\gamma^{2}}
\le20\frac{\gamma^2}{B_\gamma+
\sqrt{\gamma(1-\gamma)(r-1)}}
\le20\Phi(\gamma,r).
\]

Now, we assume that $x\ge3$, and choose $k=\lfloor x\rfloor$ in harmonic construction.  Since
$\lambda\triangleq\frac{\gamma}{1-\gamma}\le1$ and $x\ge 3$, this choice satisfies $1\le k\le r/2$.  Moreover,
\[
\frac{k}{r-k}\le\frac{x}{r-x}=\frac{x}{x^{2}/\lambda-x}=\frac{\lambda}{x-\lambda},
\qquad
\frac{\lambda}{k+\lambda}
\le\frac{\lambda}{x-1+\lambda},
\] where the first term comes from $x\ge k$ and  $r=x^{2}/\lambda$.
Eq.\eqref{eq:harmonic-upper} consequently gives
\begin{equation}\label{eq:middle-harmonic-upper}
R_{\gr}(\gamma,r)
\le\min\left\{\frac12,
2\lambda\left(
\frac1{x-\lambda}+\frac1{x-1+\lambda}
\right)\right\}.
\end{equation}
Similarly, in this setting, if $B_\gamma+
(1-\gamma)x\le 40\gamma^{2}$, we also can show that
\[
B_\gamma+
\sqrt{\gamma(1-\gamma)(r-1)}=B_\gamma+
(1-\gamma)\cdot \sqrt{\frac{\gamma}{1-\gamma}(r-1)}=B_\gamma+
(1-\gamma)x\le 40\gamma^{2},
\]
so $R_{\gr}(\gamma,r)\le\frac12=20\frac{\gamma^{2}}{40\gamma^{2}}
\le20\frac{\gamma^2}{B_\gamma+
\sqrt{\gamma(1-\gamma)(r-1)}}
\le20\Phi(\gamma,r)$.

Instead, if  $B_\gamma+
(1-\gamma)x> 40\gamma^{2}$, we define $h(x)\triangleq(B_\gamma+(1-\gamma) x)
\left(\frac1{x-\lambda}+\frac1{x-1+\lambda}\right)$. Then, from Eq.\eqref{eq:middle-harmonic-upper}, we also can have\begin{equation*}
R_{\gr}(\gamma,r)\le2\lambda\left(
\frac1{x-\lambda}+\frac1{x-1+\lambda}\right)\le2\frac{\lambda}{B_\gamma+
(1-\gamma)x}\cdot h(x).
\end{equation*}Moreover, direct differentiation of $h(x)$ gives
\[
h'(x)
=-\left(\frac{B_\gamma+(1-\gamma)\lambda}{(x-\lambda)^{2}}+\frac{B_\gamma+(1-\gamma)(1-\lambda)}{(x-1+\lambda)^{2}}\right)<0.
\]
Hence, the function $h(x)$ is decreasing and the condition $B_\gamma+
(1-\gamma)x> 40\gamma^{2}$ implies that $h(x)\le h(\frac{40\gamma^{2}-B_\gamma}{1-\gamma})$.  Moreover, due to $\lambda=\gamma/(1-\gamma)$, we also can show  
\begin{equation}\label{add_appendix_10.14}
   h(\frac{40\gamma^{2}-B_\gamma}{1-\gamma})=40\gamma^2(1-\gamma)\left(
\frac1{39\gamma^2-2}
+\frac1{39\gamma^2+3\gamma-3}
\right). 
\end{equation}
With the help of Eq.\eqref{add_appendix_10.14}, we next prove 
\begin{equation*}
  \frac1{39\gamma^2-2}
+\frac1{39\gamma^2+3\gamma-3}\le\frac{1}{4\gamma}.  
\end{equation*}
More specifically, we consider the following polynomial function $P(\gamma)$, namely,
\begin{align*}
P(\gamma)\triangleq (39\gamma^2-2)(39\gamma^2+3\gamma-3)
-4\gamma(78\gamma^2+3\gamma-5)=
1521\gamma^4-195\gamma^3-207\gamma^2+14\gamma+6.
\end{align*}
Then, we also can show that, when $\frac1{2\sqrt2}\le\gamma\le\frac12$, the following identities hold:
\[
P'''(\gamma)=36504\gamma-1170>0,
\quad
P''\!\left(\frac1{2\sqrt2}\right)
=\frac{3735-585\sqrt2}{2}>0,
\]
and
\[
P'\!\left(\frac1{2\sqrt2}\right)
=\frac{693\sqrt2-473}{8}>0.
\]
Thus $P$ is increasing over $\frac1{2\sqrt2}\le\gamma\le\frac12$, and $P(\frac{1}{2\sqrt{2}})=(249-166\sqrt2)/64>0$. Thus, we get  $ \frac1{39\gamma^2-2}
+\frac1{39\gamma^2+3\gamma-3}\le\frac{1}{4\gamma}$ such that\begin{align*}
R_{\gr}(\gamma,r)&\le2\lambda\left(
\frac1{x-\lambda}+\frac1{x-1+\lambda}\right)\le2\frac{\lambda}{B_\gamma+
(1-\gamma)x}\cdot h(x)\le2\frac{\lambda}{B_\gamma+
(1-\gamma)x}\cdot h(\frac{40\gamma^{2}-B_\gamma}{1-\gamma})\\&\le 20\frac{\gamma^{2}}{B_\gamma+
(1-\gamma)x}\le 20\Phi(\gamma,r).
\end{align*}

Finally, we consider that $1/2\le\gamma<1$. In this setting, $\lambda=\gamma/(1-\gamma)\ge1$ such that Eq.\eqref{add_eq:harmonic-fixed-rate} becomes
\begin{equation}\label{eq:near-upper-harmonic}
R_{\gr}(\gamma,r)
\le\frac1{\sqrt r}
\max\left\{
8\sqrt{\frac{\gamma}{1-\gamma}},
2\sqrt{\frac{1-\gamma}{\gamma}}
\right\}.  
\le8\sqrt{\frac{\gamma}{(1-\gamma)r}}=\frac{8\gamma}{\sqrt{\gamma(1-\gamma)r}}.
\end{equation}
If $\sqrt{\gamma(1-\gamma)r}\ge16\gamma$, then $\sqrt{\gamma(1-\gamma)r}>8$, so
$B_\gamma\le2<\sqrt{\gamma(1-\gamma)r}/4$ and $\sqrt{\gamma(1-\gamma)(r-1)}\le\sqrt{(1-\gamma)r}$. Then, from
Eq.\eqref{eq:near-upper-harmonic}, we can show that
\begin{align*}
R_{\gr}(\gamma,r)&\le\frac{8\gamma}{\sqrt{\gamma(1-\gamma)r}}=\frac{16\cdot(\gamma/2)}{\sqrt{\gamma(1-\gamma)r}}\le\frac{16\gamma^2}{\sqrt{\gamma(1-\gamma)r}}=\frac{20\gamma^2}{\sqrt{\gamma(1-\gamma)r}/4+\sqrt{\gamma(1-\gamma)r}}
\\&\le20\frac{\gamma^2}{B_\gamma+\sqrt{\gamma(1-\gamma)(r-1)}}\le20\Phi(\gamma,r).
\end{align*}

When $\sqrt{\gamma(1-\gamma)r}\le16\gamma$, then 
\[
B_\gamma+
\sqrt{\gamma(1-\gamma)(r-1)}\le2+16\gamma\le40\gamma^2,
\]
where the last inequality comes from that $40\gamma^2-16\gamma-2\ge0$ for every $\gamma\ge1/2$.  Consequently,
\[
R_{\gr}(\gamma,r)
\le\frac12
\le20\frac{\gamma^2}{B_\gamma+
\sqrt{\gamma(1-\gamma)(r-1)}}
\le20\Phi(\gamma,r).
\]
Finally, we have $R_{\gr}(\gamma,r)\le20\Phi(\gamma,r)$ for every
$0<\gamma<1$.  At $\gamma=1$, the classical tight submodular guarantee
gives $R_{\gr}(1,r)=1/2$.  Thus, we prove Eq.\eqref{eq:uniform-upper}.
\end{proof}

\section{Proofs for the lower bounds}
\label{app:lower-proofs}
\addappendixsection{Proofs for the lower bounds}
\subsection{Proof of Lemma~\ref{lem:basic-identities}}
\label{app:proof-basic-identities}
\addappendixsubsection{Proof of Lemma~\ref{lem:basic-identities}}
\begin{proof}
Since $S_{t-1}\subseteq S_t$, the monotonicity of $f$ gives
$d_t=f(S_t)-f(S_{t-1})\ge0$, proving the first claim.  For the second one, according to
$d_i=F_i-F_{i-1}$, we have 
\begin{align*}
W_t=\sum_{i=1}^t i(F_i-F_{i-1})=\sum_{i=1}^t iF_i-\sum_{i=1}^t iF_{i-1}
=\sum_{i=1}^t iF_i-\sum_{j=0}^{t-1}(j+1)F_j=tF_t-\sum_{j=0}^{t-1}F_j
=tF_t-K_t.
\end{align*}
Finally, since $i\ge1$ and $d_i\ge0$ for every $i$, $W_t=\sum_{i=1}^t i d_i
\ge\sum_{i=1}^t d_i
=F_t-F_0=F_t.$
\end{proof}


\subsection{Proof of Lemma~\ref{lem:delayed}}
\label{app:proof-delayed}
\addappendixsubsection{Proof of Lemma~\ref{lem:delayed}}
\begin{proof}
At first, we assume that the optimal base can be ordered as $O=\{o_1,\ldots,o_r\}$ and satisfies $S_t\cup\{o_{t+1},\ldots,o_r\}$
for any $t=0,\dots,r-1$. Then, we set $A=S_{i-1}$ and $C=S_t\cup\{o_i\}$. The new elements in $C\setminus A$ belong to
\(\{s_i,\ldots,s_t,o_i\}\), which contains at most $t-i+2$ elements.
Lemma~\ref{lem:ordering} shows that $S_{i-1}+o_i\in\cI$, so $o_i$ is
feasible at $S_{i-1}$.  Moreover, for every $j\ge i$,
$S_{i-1}+s_j\subseteq S_j\in\cI$.  The downward-closed property of a
matroid therefore gives $S_{i-1}+s_j\in\cI$, so $s_j$ is also feasible at
$S_{i-1}$. The standard forward greedy algorithm selects $s_i$ at iteration
$i$, and hence every such feasible element has marginal value at most
$d_i$.  Consequently,
\[
\sum_{e\in C\setminus A}f(e\mid S_{i-1})
\le(t-i+2)d_i.
\]
On the other hand, weak submodularity gives
\[
\gamma\bigl(f(C)-f(A)\bigr)
\le\sum_{e\in C\setminus A}f(e\mid S_{i-1}),
\]
and the joint gain has the exact decomposition
\[
f(C)-f(A)=F_t-F_{i-1}+f(o_i\mid S_t).
\]
Combining the last three displays gives
\[
\gamma\bigl(F_t-F_{i-1}+f(o_i\mid S_t)\bigr)
\le \sum_{e\in C\setminus A}f(e\mid S_{i-1})
\le (t-i+2)d_i.
\]
After subtracting $\gamma(F_t-F_{i-1})$, we prove
Eq.\eqref{eq:delayed}.
\end{proof}

\subsection{Proof of Lemma~\ref{lem:prefix}}
\label{app:proof-prefix}
\addappendixsubsection{Proof of Lemma~\ref{lem:prefix}}
\begin{proof}
We first sum the results of Lemma~\ref{lem:delayed}, namely,
\begin{align}
\gamma\sum_{i=1}^t f(o_i\mid S_t)
&\le \sum_{i=1}^t
\bigl((t-i+2)d_i-\gamma(F_t-F_{i-1})\bigr) \notag\\
&=(t+2)\sum_{i=1}^t d_i-\sum_{i=1}^t i d_i
  -\gamma\left(tF_t-\sum_{i=1}^tF_{i-1}\right) \notag\\
&=(t+2)F_t-W_t-\gamma(tF_t-K_t)=(t+2)F_t-(1+\gamma)W_t, \label{eq:blocked-sum}
\end{align} where the final equality follows from $tF_t-K_t=W_t$ of Lemma~\ref{lem:basic-identities}.

For every $i>t$, Lemma~\ref{lem:ordering} shows that $o_i$ is feasible at
$S_t$. The element chosen by the standard forward greedy algorithm at
iteration $t+1$ has marginal value $d_{t+1}$.  Hence, for $t<r$,
\begin{equation}\label{eq:unblocked-sum}
\sum_{i=t+1}^r f(o_i\mid S_t)\le(r-t)d_{t+1}.
\end{equation}
Moreover, throughout weak submodularity, we have 
\[
\sum_{i=1}^r f(o_i\mid S_t)\ge \gamma\bigl(f(S_t\cup O)-F_t\bigr)\ge \gamma(\OPT-F_t).
\]
Therefore, for $1\le t\le r$, Eq.\eqref{eq:blocked-sum} and
Eq.\eqref{eq:unblocked-sum} imply
\begin{align*}
\gamma^2(\OPT-F_t)
&\le \gamma\sum_{i=1}^t f(o_i\mid S_t)
   +\gamma\sum_{i=t+1}^r f(o_i\mid S_t)\\
&\le (t+2)F_t-(1+\gamma)W_t
   +\gamma(r-t)d_{t+1}.
\end{align*}
After adding $\gamma^2F_t$, we can obtain
\begin{align*}
\gamma^2\OPT
&\le (t+2+\gamma^2)F_t-(1+\gamma)W_t
   +\gamma(r-t)d_{t+1}\\
&=(2-\gamma+\gamma^2)F_{t}+\left(tF_{t}-W_{t}\right)+\gamma(F_{t}-W_{t})+\gamma(r-t)d_{t+1}
 \\  
&=B_\gamma F_t+K_t-\gamma(W_t-F_t)
   +\gamma(r-t)d_{t+1}\le B_\gamma F_t+K_t+\gamma(r-t)d_{t+1},
\end{align*}where the final inequality follows from $W_t\ge F_t$ of Lemma~\ref{lem:basic-identities}. Thus, we prove Eq.\eqref{eq:prefix-main} for $1\le t\le r$.  

Particularly when $t=r$, 
\begin{align*}
\gamma^2\OPT\le B_\gamma F_r+K_r,
\end{align*}
which proves Eq.\eqref{eq:prefix-terminal}.  Finally, at $t=0$, the weak
submodularity implies that
\[
\gamma\OPT
\le\sum_{o\in O}f(o\mid\varnothing)
\le r d_1.
\]
Multiplying by $\gamma$ proves Eq.\eqref{eq:prefix-main} for $t=0$. Note that $f(\emptyset)=0$.
\end{proof}
\subsection{Proof of Lemma~\ref{lem:recurrence}}
\label{app:proof-recurrence}
\addappendixsubsection{Proof of Lemma~\ref{lem:recurrence}}
\begin{proof}
Set $V=C-BM$.  If $V\le0$, then $C\le BM\le\left(B+\sqrt{A(1-1/r)}\right)M$. So we suppose $V>0$. From $C\le BM+K_r$, we know that $K_r\ge V$.  Let $h$ be the least index such that
$K_h\ge V$.  For $t<h$,
\[
V-K_t
=C-BM-K_t
\le A(F_{t+1}-F_t)+B(F_t-M)
\le A(F_{t+1}-F_t).
\]
Since $K_t<V$ for $t<h$,
both sides of the previous inequality are nonnegative.  Multiplying by $F_t\ge0$ and summing gives
\begin{equation}\label{eq:recurrence-summed}
\sum_{t=0}^{h-1}(V-K_t)F_t
\le A\sum_{t=0}^{h-1}F_t(F_{t+1}-F_t).
\end{equation}
Moreover, due to $K_{t+1}=K_t+F_t$, we can show $K_{t+1}^2-K_t^2=2K_tF_t+F_t^2$. Summing this identity yields
\[
\sum_{t<h}K_tF_t
=\frac12\left(K_h^2-\sum_{t<h}F_t^2\right).
\]
Consequently,
\begin{align*}
\sum_{t<h}(V-K_t)F_t=VK_h-\frac12K_h^2+
\frac12\sum_{t<h}F_t^2=\frac12V^2-\frac12(K_h-V)^2+
\frac12\sum_{t<h}F_t^2.
\end{align*}
The minimality of $h$ gives
$K_{h-1}<V\le K_h=K_{h-1}+F_{h-1}$, and hence
$0\le K_h-V\le F_{h-1}$ such that  $(K_h-V)^2\le F_{h-1}^2\le\sum_{t<h}F_t^2$. So we have 
\begin{equation}\label{eq:recurrence-left}
\sum_{t<h}(V-K_t)F_t=\frac12V^2+
\frac12\left(\sum_{t<h}F_t^2-(K_h-V)^2\right)\ge\frac12V^2.
\end{equation}

Next, we investigate the right-hand side of Eq.\eqref{eq:recurrence-summed}. Due to $d_{t+1}\triangleq F_{t+1}-F_{t}$, we can show that $2F_td_{t+1}=F_{t+1}^2-F_t^2-d_{t+1}^2$. Summing this identity for $t=0,\ldots,h-1$ and using $F_0=0$, we then have
\begin{align*}
A\sum_{t<h}F_t(F_{t+1}-F_t)
&=A\sum_{t<h}F_td_{t+1}=\frac A2\left(
\sum_{t<h}(F_{t+1}^2-F_t^2)-\sum_{t<h}d_{t+1}^2
\right)\\&=\frac A2\left(F_h^2-\sum_{t<h}d_{t+1}^2\right)
\le\frac A2\left(1-\frac1h\right)F_h^2
\le\frac A2\left(1-\frac1r\right)M^2,
\end{align*} where the first inequality follows from the Cauchy-Schwarz inequality $\sum_{t<h}d_{t+1}^2
\ge\frac1h\left(\sum_{t<h}d_{t+1}\right)^2
=\frac{F_h^2}{h}$ and the second inequality uses $F_h\le M$.  Thus
\begin{equation}\label{eq:recurrence-right}
A\sum_{t<h}F_t(F_{t+1}-F_t)
\le\frac A2\left(1-\frac1r\right)M^2.
\end{equation}
Merging Eq.\eqref{eq:recurrence-left} and
Eq.\eqref{eq:recurrence-right} into Eq.\eqref{eq:recurrence-summed}, we obtain
\[
\frac12V^2
\le\sum_{t<h}(V-K_t)F_t
\le A\sum_{t<h}F_t(F_{t+1}-F_t)
\le\frac A2\left(1-\frac1r\right)M^2.
\]
Therefore $V=C-BM\le\sqrt{A(1-1/r)}M$, which proves
Eq.\eqref{eq:recurrence-conclusion}.
\end{proof}

\subsection{Proof of the rank-dependent terms in Theorem~\ref{thm:lower}.}
\label{app:proof-universal}
\addappendixsubsection{Proof of  Theorem~\ref{thm:lower}}
\begin{proof}
In this section, we show the verification on Theorem~\ref{thm:lower}. At first, we set 
\[
A\triangleq \gamma r,\qquad B\triangleq B_\gamma,
\qquad C\triangleq \gamma^2\OPT,\qquad M\triangleq F_r.
\]
Lemma~\ref{lem:basic-identities} shows that
$0=F_0\le F_1\le\cdots\le F_r=M$.  Moreover, for every $t<r$,
Lemma~\ref{lem:prefix} gives
\begin{align*}
C\le B_\gamma F_t+K_t+\gamma(r-t)(F_{t+1}-F_t)\le BF_t+K_t+A(F_{t+1}-F_t),
\end{align*}
where the second inequality uses
$\gamma(r-t)\le A$.  The same lemma also gives $C\le BM+K_r$ through
Eq.\eqref{eq:prefix-terminal}.  Thus, by Lemma~\ref{lem:recurrence}, we have
yields
\[
\gamma^2\OPT
\le\left(B_\gamma+\sqrt{\gamma r(1-1/r)}\right)F_r
=\left(B_\gamma+\sqrt{\gamma(r-1)}\right)F_r.
\]
Consequently, $\frac{F_r}{\OPT}
\ge\frac{\gamma^2}{B_\gamma+\sqrt{\gamma(r-1)}}$. 

Moreover, every element of $O$ is feasible at the empty set, and
the first standard forward greedy selection has marginal value $F_1$. Therefore,
\[
rF_1\ge\sum_{o\in O}f(o\mid\varnothing)
\ge\gamma\OPT,
\]
where the second inequality follows from weak submodularity applied to
$\varnothing\subseteq O$.  Monotonicity gives
$\frac{F_r}{\OPT}\ge \frac{F_1}{\OPT}\ge\gamma/r$.  Combining these two bounds, we prove Eq.\eqref{eq:main-lower}.
\end{proof}

\subsection{Proof of Theorem~\ref{thm:endpoint}}
\label{app:proof-endpoint}
\addappendixsubsection{Proof of Theorem~\ref{thm:endpoint}}
\begin{proof} Firstly, for any set $A$ and two
distinct elements $s,e\notin A$, $\gamma$-weak submodularity implies that
\[
f(s\mid A)+f(e\mid A)
\ge\gamma\bigl(f(A+s+e)-f(A)\bigr)=\gamma\bigl(f(s\mid A)+f(e\mid A+s)\bigr).
\]
Then, if $\gamma>0$, we can obtain
\begin{equation}\label{eq:pairwise-growth}
f(e\mid A+s)
\le\gamma^{-1}f(e\mid A)
+\frac{1-\gamma}{\gamma}f(s\mid A).
\end{equation}

Let $O=\{o_1,\ldots,o_r\}$ be ordered as in
Lemma~\ref{lem:ordering}.  We next bound $f(o_i\mid S_r)$ for every
$i\in[r]$.  If $o_i\in S_r$, then $f(o_i\mid S_r)=0$. Therefore, we suppose that $o_i\notin S_r$, and define
\[
x_t\triangleq f(o_i\mid S_t),
\qquad t=i-1,i,\ldots,r.
\]
Lemma~\ref{lem:ordering} shows that $o_i$ is feasible at $S_{i-1}$.
The greedy choice at iteration $i$ therefore implies
\begin{equation}\label{eq:endpoint-initial}
x_{i-1}=f(o_i\mid S_{i-1})\le d_i.
\end{equation}
Then, for each $t=i,\ldots,r$, we apply Eq.\eqref{eq:pairwise-growth} with
$A=S_{t-1}$, $s=s_t$, and $e=o_i$.  Since
$S_t=S_{t-1}+s_t$ and $f(s_t\mid S_{t-1})=d_t$, this gives the following inequality:
\begin{equation}\label{eq:endpoint-recurrence}
x_t
\le\gamma^{-1}x_{t-1}
+\frac{1-\gamma}{\gamma}d_t.
\end{equation}

We now iterate Eq.\eqref{eq:endpoint-recurrence}. We can get, for any $m=i,\ldots,r$,
\begin{equation}\label{eq:endpoint-iteration}
x_m
\le\gamma^{-(m-i+1)}d_i
+(1-\gamma)\sum_{t=i}^m
\gamma^{-(m-t+1)}d_t.
\end{equation}
Taking $m=r$, we can show that, for every $i$,
\begin{equation}\label{eq:endpoint-terminal-marginal}
f(o_i\mid S_r)
\le\gamma^{-(r-i+1)}d_i
+(1-\gamma)\sum_{t=i}^r
\gamma^{-(r-t+1)}d_t.
\end{equation}

We sum Eq.\eqref{eq:endpoint-terminal-marginal} over $i=1,\ldots,r$.
For a fixed $t$, the first term contributes
$\gamma^{-(r-t+1)}d_t$ once.  In the second term of Eq.\eqref{eq:endpoint-terminal-marginal}, $d_t$ only occurs for
exactly the indices $i=1,\ldots,t$, and therefore contributes
$t(1-\gamma)\gamma^{-(r-t+1)}d_t$.  Hence, we have 
\begin{equation}\label{eq:endpoint-summed}
\sum_{i=1}^r f(o_i\mid S_r)
\le\sum_{t=1}^r
\gamma^{-(r-t+1)}\bigl[1+t(1-\gamma)\bigr]d_t.
\end{equation}

It remains to bound the coefficient of each $d_t$. Firstly, we set $a_t=\gamma^{-(r-t+1)}\bigl(1+t(1-\gamma)\bigr)$.Then, we show that 
$a_t$ is nonincreasing because
\[
\frac{a_{t+1}}{a_t}
=1-\frac{(t+1)(1-\gamma)^2}{1+t(1-\gamma)}
\le1.
\]
Thus $a_t\le a_1=\gamma^{-r}(2-\gamma)$, and we have 
\begin{align}
\sum_{i=1}^r f(o_i\mid S_r)\le(2-\gamma)\gamma^{-r}\sum_{t=1}^r d_t=(2-\gamma)\gamma^{-r}F_r.\label{eq:endpoint-upper-sum}
\end{align}

Finally, $\gamma$-weak submodularity gives
\begin{align*}
\gamma(\OPT-F_r)\le\gamma\bigl(f(S_r\cup O)-F_r\bigr)\le\sum_{o\in O\setminus S_r}f(o\mid S_r)
=\sum_{i=1}^r f(o_i\mid S_r).
\end{align*}
Combining this inequality with Eq.\eqref{eq:endpoint-upper-sum} yields $\gamma(\OPT-F_r)
\le(2-\gamma)\gamma^{-r}F_r$. Thus, we obtain
\[
\frac{F_{r}}{\OPT}\ge\frac{\gamma}{(2-\gamma)\gamma^{-r}+\gamma}=\frac{\gamma^{r+1}}{2-\gamma+\gamma^{r+1}},
\]
which proves Eq.\eqref{eq:endpoint-lower}.
\end{proof}

\subsection{Proof of Lemma~\ref{lem:lower-comparison}}
\label{app:lower-comparison}
\addappendixsubsection{Proof of Lemma~\ref{lem:lower-comparison}}
\begin{proof}
Set $x\triangleq\sqrt{\gamma r}$.  Since
$B_\gamma=2-\gamma+\gamma^2\le2$ and
$\sqrt{\gamma(r-1)}<\sqrt{\gamma r}=x$, we have
\[
B_\gamma+\sqrt{\gamma(r-1)}
<2+x
<2+2x+\gamma^2.
\]
This proves the first strict inequality in
Eq.\eqref{eq:lower-comparison}.

For the second one, it is worth noting that, when $r\ge2$, we have 
$x\ge\sqrt{2\gamma}>\gamma$.  Therefore,
\[
\frac12(1+x)>\frac12(1+\gamma)\ge\gamma\ge\gamma^2,
\]
and hence
\[
\frac{\gamma^{2}}{2+2x+\gamma^2}
>\frac{\gamma^{2}}{2+2x+\frac12(1+x)}
=\frac{\gamma^{2}}{\frac{5}{2}(1+x)}=\frac{0.4\gamma^{2}}{1+x}.
\]
\end{proof}

\subsection{Proof of Lemma~\ref{lem:near-prefix}}
\label{app:proof-near-prefix}
\addappendixsubsection{Proof of Lemma~\ref{lem:near-prefix}}
In this section, we prove Lemma~\ref{lem:near-prefix}. Before that, we first give a refined result analogous to Lemma~\ref{lem:delayed}, that is to say,

\begin{lemma}
\label{lem:near-delayed}
For any $1/2\le\gamma<1$ and $1\le i\le t\le r$, we have 
\begin{equation}\label{eq:near-delayed}
\gamma f(o_i\mid S_t)
\le15\bigl(3+(1-\gamma)(t-i+1)\bigr)d_i.
\end{equation}
\end{lemma}

With the help of Eq.\eqref{eq:near-delayed}, we then verify Lemma~\ref{lem:near-prefix}, that is to say,
\begin{proof}For any $1\le t\le r$,  summing
Eq.\eqref{eq:near-delayed} over $i=1,\ldots,t$ gives
\begin{align*}
\gamma\sum_{i=1}^t f(o_i\mid S_t)
&\le45\sum_{i=1}^t d_i
+15(1-\gamma)\sum_{i=1}^t(t-i+1)d_i\\
&=45\sum_{i=1}^t d_i
+15(1-\gamma)\left((t+1)\sum_{i=1}^td_i-\sum_{i=1}^tid_i\right)\\
&=45F_t
+15(1-\gamma)\left((t+1)F_t-W_t\right)\\&=45F_t+15(1-\gamma)(F_t+K_t)=(45+15(1-\gamma))F_t+15(1-\gamma) K_t,
\end{align*} where the third equality comes from $W_t=tF_t-K_t$ in Lemma~\ref{lem:basic-identities}.

For every $i>t$, Lemma~\ref{lem:ordering} shows that $o_i$ is feasible at
$S_t$. The element chosen by the standard forward greedy algorithm at
iteration $t+1$ has marginal value $d_{t+1}$.  Hence, for $t<r$,
\begin{equation*}
\sum_{i=t+1}^r f(o_i\mid S_t)\le(r-t)d_{t+1}.
\end{equation*}As a result, we have 
\begin{align*}
\gamma^2(\OPT-F_t)\le\gamma\sum_{i=1}^r f(o_i\mid S_t)\le(45+15(1-\gamma))F_t+15(1-\gamma)K_t
+\gamma(r-t)d_{t+1}.
\end{align*}
After adding $\gamma^2F_t$ to both sides, we get
\begin{equation*}
    \begin{aligned}
        \gamma^{2}\OPT&\le(\gamma^{2}+45+15(1-\gamma))F_t+15(1-\gamma)K_t
+\gamma(r-t)d_{t+1}\\&=D_\gamma F_t+15(1-\gamma)K_t
+\gamma(r-t)d_{t+1},
    \end{aligned}
\end{equation*} where  $\gamma^2+45+15(1-\gamma)
=60-15\gamma+\gamma^2
=D_\gamma$.

Particularly when $t=r$,  we have 
\[
\gamma^2\OPT
\le D_\gamma F_r+15(1-\gamma)K_r,
\]
which is Eq.\eqref{eq:near-prefix-terminal}. 
\end{proof}

In the remainder of this section, we show the proof of Lemma~\ref{lem:near-delayed}.
\subsubsection{Proof of Lemma~\ref{lem:near-delayed}}
\label{app:proof-near-delayed}
\begin{proof}
Fix $1\le i\le t\le r$.  If $o_i\in S_t$, then
$f(o_i\mid S_t)=0$, and then the claim of Lemma~\ref{lem:near-delayed} immediately follows. Thus, we suppose 
that $o_i\notin S_t$, and set
\[
h\triangleq t-i+1,\qquad
x_\ell\triangleq f(o_i\mid S_{i-1+\ell})\quad(0\le\ell\le h).
\]
Moreover, Lemma~\ref{lem:ordering} shows that $o_i$ is feasible at $S_{i-1}$, so the
greedy choice at iteration $i$ gives
\begin{equation}\label{eq:near-initial}
x_0\le d_{i}.
\end{equation}

Next, for any $0\le u<v\le h$, we apply $\gamma$-weak submodularity with $A=S_{i-1+u}$ and $C=S_{i-1+v}\cup\{o_i\}$, that is to say,
\begin{equation*}
    \sum_{e\in C\setminus A}f(e\mid A)\ge\gamma\left(f(C)-f(A)\right).
\end{equation*}Importantly, we know that, $C\setminus A=\{s_{i+u},\dots,s_{i+v-1}\}\cup\{o_i\}$. Moreover, for every $j\in\{u+1,\ldots,v\}$, the element
$s_{i+j-1}$ is feasible at $A$, because $A+s_{i+j-1}\subseteq S_{i+j-1}\in\cI$.Therefore, we have $f(s_{i+j-1}\mid A)\le d_{i+u}$ for any $u+1\le j\le v$. As a result, we can show
\[
\sum_{e\in C\setminus A}f(e\mid A)=f(o_{i}\mid S_{i-1+u})+\sum_{j=u+1}^{v}f(s_{i+j-1}\mid S_{i-1+u})\le x_u+(v-u)d_{i+u}.
\] Similarly, the joint gain $f(C)-f(A)$ can be represented as follows:
\[
f(C)-f(A)=f(o_{i}\mid S_{i-1+v})+\sum_{j=u+1}^{v}d_{i+j-1}=x_v+\sum_{j=u+1}^{v}d_{i+j-1}.
\]
Hence
\[
x_u+(v-u)d_{i+u}
\ge\gamma\left(x_v+\sum_{j=u+1}^{v}d_{i+j-1}\right),
\]
or equivalently,
\begin{equation}\label{eq:near-batch}
\gamma x_v-x_u
+\gamma\sum_{j=u+1}^{v}d_{i+j-1}
-(v-u)d_{i+u}\le0.
\end{equation}

In the remainder of this section, we will construct a nonnegative linear combination of instances of Eq.\eqref{eq:near-batch}, indexed by suitable choices of $u$ and $v$, which can cause the coefficients of all intermediate marginals $x_1,\ldots,x_{h-1}$ to vanish while ensuring that those of the later greedy increments $d_{i+1},\ldots,d_{i+h-1}$ are nonnegative. After discarding these nonnegative terms, the resulting inequality involves
only $x_0$, $x_h$, and $d_i$.

\paragraph{Weight construction and coefficient cancellation.} Before going into the details, we first define some useful symbols, namely,
\[
M\triangleq 3+(1-\gamma)h,
\qquad z\triangleq 1+\frac{1-\gamma}{M},
\qquad \rho\triangleq z^{-1},
\qquad L\triangleq\lceil5M\rceil,\qquad G\triangleq\sum_{a=0}^{L-1}z^{a}.
\] With these notations, we next show two frequently used results, i.e., $\gamma G\ge L$ and $\gamma z^L\ge1$. More specifically, 
\begin{equation*}
\begin{aligned}
G&=\sum_{a=0}^{L-1}\left(1+\frac{1-\gamma}{M}\right)^a\ge
\sum_{a=0}^{L-1}\left(1+\frac{a(1-\gamma)}{M}\right)=
L\left(1+\frac{(1-\gamma)(L-1)}{2M}\right)\\&
\ge L\left(1+\frac{(1-\gamma)(5M-1)}{2M}\right)=L\left(1+(1-\gamma)\left(\frac{5}{2}-\frac{1}{2M}\right)\right)\ge L\left(1+(1-\gamma)\frac{1}{\gamma}\right)=\frac{L}{\gamma},
\end{aligned}
\end{equation*} where the first inequality follows from $(1+q)^{a}\ge 1+aq$ and the final inequality comes from $M\ge 3$ and $\gamma\ge 1/2$ such that $\frac{5}{2}-\frac{1}{2M}\ge\frac{7}{3}\ge2\ge\frac{1}{\gamma}$.

After that, due to $G=\frac{1-z^{L}}{1-z}$, we can show that 
\begin{equation*}
   z^L=1+(z-1)G
=1+\frac{1-\gamma}{M}G
\ge1+\frac{(1-\gamma)L}{\gamma M}
\ge1+\frac{1-\gamma}{\gamma}
=\frac1\gamma, 
\end{equation*}where the first inequality follows from the previous $G\ge\frac{L}{\gamma}$  and the final one comes from $L\ge M$. 

Subsequently, for any $k\in\{1-L,\dots,h-1\}$, we consider the following two endpoints
\[u_k\triangleq\max\{0,k\},\qquad
v_k\triangleq\min\{h,k+L\},\] and assign each equation $B_{k}\triangleq \gamma x_{v_k}-x_{u_k}
+\gamma\sum_{j=u_k+1}^{v_k}d_{i+j-1}
-(v_k-u_k)d_{i+u_k}$  with the weight $w_{k}\triangleq\frac{z^{h-1-k}}{G}$. As a result, we can easily show
\begin{equation}\label{equ_add_10_15}\sum_{k=1-L}^{h-1}w_{k}B_{k}-\sum_{j=1}^{h}(\gamma z^{L}-1)w_jx_j\le 0,
\end{equation} where $B_k\le 0$ from Eq.\eqref{eq:near-batch} and $\gamma z^{L}\ge1$. 

Then, we compute the coefficients of $x_0,\ldots,x_h$ in Eq.\eqref{equ_add_10_15}. For any $1\le j<h$, the term $x_j$ only appears in $B_{j-L}$, with coefficient $\gamma w_{j-L}$, and in $B_j$, with coefficient $-w_j$. In addition, the
auxiliary term $-\sum_{j=1}^{h}(\gamma z^{L}-1)w_jx_j$ contributes the coefficient $-(\gamma z^{L}-1)w_j$. Thus, the total coefficient of $x_j$ is
\begin{align*}
\gamma w_{j-L}-w_j-(\gamma z^L-1)w_j
=\gamma z^Lw_j-w_j-(\gamma z^L-1)w_j=0,
\end{align*} where the identity follows from $w_{j-L}=z^Lw_j$. So all $x_1,\ldots,x_{h-1}$ will be canceled exactly.

It remains to compute the coefficients of the two boundary marginals. The term $x_h$ appears in $B_k$ with coefficient $\gamma w_k$ for any
$k\in\{h-L,\ldots,h-1\}$. Thus, the total coefficient of $x_h$ is
\[
\gamma\sum_{k=h-L}^{h-1}w_k
=
\gamma\frac{1}{G}\sum_{a=0}^{L-1}z^a
=\gamma.
\]
Similarly, $x_0$ appears in $B_k$ with coefficient $- w_k$ for any
$k\in\{1-L,\ldots,0\}$. Thus, the total coefficient of $x_0$ is
\begin{equation}\label{add_Y}
-\sum_{k=1-L}^{0}w_k
=\frac{z^{h-1}}{G}\sum_{k=1-L}^{0}z^{-k}
=
-\frac{z^{h-1}}{G}\sum_{a=0}^{L-1}z^a
=
-z^{h-1}.    
\end{equation}

Next, we compute the coefficients of $d_{i},\ldots,d_{i+h-i}$ in
Eq.~\eqref{equ_add_10_15}. Fix any $2\le j\le h$. The term $d_{i+j-i}$ appears
in the positive sum $\gamma\sum_{j=u_k+1}^{v_k}d_{i+j-1}$ of $B_k$ for $k\in\{j-L,j-L+1,\ldots,j-1\}$.
Thus, its total positive coefficient is $\gamma\sum_{k=j-L}^{j-1}w_k
=\gamma w_{j-1}\sum_{a=0}^{L-1}z^a=
\gamma w_{j-1}G$. On the other hand, the negative contribution to $d_{i+j-i}$ comes from the term
$B_{j-1}$ with coefficient $-(v_{j-1}-u_{j-1})w_{j-1}$. Since $v_{j-1}-u_{j-1}\le L$, the total coefficient $c_j$ of $d_{i+j-1}$
satisfies
\begin{equation*}
c_j=\gamma w_{j-1}G-(v_{j-1}-u_{j-1})w_{j-1}
\ge
\bigl(\gamma G-L\bigr)w_{j-1}
\ge 0,
\qquad 2\le j\le h,
\end{equation*}
where the last inequality follows from $\gamma G\ge L$.

It remains to compute the coefficient of $d_i$. The term $d_i$ appears in the positive sum $\gamma\sum_{j=u_k+1}^{v_k}d_{i+j-1}$ of $B_k$ for any $k\in\{1-L,\ldots,0\}$. Thus, its total positive coefficient is  $\gamma\sum_{k=1-L}^{0}w_k$. Similarly, the negative terms contribute
$-\sum_{k=1-L}^{0}v_kw_k$. Therefore, the total coefficient $C$ of
$d_i$ is
\[
C=\gamma\sum_{k=1-L}^{0}w_k-\sum_{k=1-L}^{0}v_kw_k\ge(\gamma-L)\sum_{k=1-L}^{0}w_k=(\gamma-L)z^{h-1},
\] where the first inequality comes from $v_k\le L$ and the final equality follow from Eq.\eqref{add_Y}.

Combining the preceding coefficient calculations, we obtain
\[
\sum_{k=1-L}^{h-1}w_{k}B_{k}-\sum_{j=1}^{h-1}(\gamma z^{L}-1)w_jx_j=\gamma x_h-z^{h-1}\cdot x_0+C\cdot d_i
+\sum_{j=2}^{h-1}c_jd_{i+j-1}
\le 0.
\]Due to $c_j\ge0$, we have 
\begin{align*}
\gamma x_h\le z^{h-1}\cdot x_0-C\cdot d_i\le
(z^{h-1}-C)d_i\le
(1-\gamma+L)z^{h-1}\cdot d_i.
\end{align*}

Finally, $z^{h-1}
=
\left(1+\frac{1-\gamma}{M}\right)^{h-1}
\le
\exp\left(\frac{(1-\gamma)(h-1)}{M}\right)
\le e$. Moreover, since $L=\lceil5M\rceil$, $\gamma\ge1/2$, and $M\ge3$, $1-\gamma+L
\le
\frac{3}{2}+5M
\le
\frac{11}{2}M$. Consequently, we have  $\gamma x_{h}\le(1-\gamma+L)z^{h-1}\cdot d_i\le\frac{11e}{2}M\cdot d_i\le 15M\cdot d_i $. Recalling that $M=3+(1-\gamma)h$ and $h=t-i+1$, we thus conclude
that
\[
\gamma f(o_i\mid S_t)=\gamma x_h\le15M\cdot d_i =
15\bigl(3+(1-\gamma)(t-i+1)\bigr)d_i.
\]
\end{proof}

\subsection{Proof of Theorem~\ref{thm:near-lower}}
\label{app:proof-near-lower}
\addappendixsubsection{Proof of Theorem~\ref{thm:near-lower}}
\begin{proof}
The term $\gamma/r$ follows from Theorem~\ref{thm:lower}.  We prove the
second term.  We set 
\[
H\triangleq 15(1-\gamma),
\qquad \widetilde F_t\triangleq HF_t,
\qquad \widetilde K_t\triangleq HK_t,
\qquad \widetilde M\triangleq HF_r.
\]
Because $1/2\le\gamma<1$, we have $H>0$.  Lemma~\ref{lem:near-prefix}
implies, for every $t<r$,
\begin{align*}
\gamma^2\OPT
\le D_\gamma F_t+HK_t+\gamma(r-t)(F_{t+1}-F_t)\le\frac{D_\gamma}{H}\widetilde F_t+\widetilde K_t
+\frac{\gamma r}{H}
(\widetilde F_{t+1}-\widetilde F_t).
\end{align*}
Particularly when $t=r$, we also have $\gamma^2\OPT
\le\frac{D_\gamma}{H}\widetilde M+\widetilde K_r$.

Then, if we apply Lemma~\ref{lem:recurrence} with
\[
A=\frac{\gamma r}{H},
\qquad B=\frac{D_\gamma}{H},
\qquad C=\gamma^2\OPT,
\qquad M=\widetilde M,
\] we can have \begin{align*}
\gamma^2\OPT
\le\left(
\frac{D_\gamma}{H}
+\sqrt{\frac{\gamma r}{H}\left(1-\frac1r\right)}
\right)\widetilde M=\left(D_\gamma+\sqrt{H\gamma(r-1)}\right)F_r=\left(D_\gamma+
\sqrt{15\gamma(1-\gamma)(r-1)}\right)F_r.
\end{align*}
Therefore, $\frac{F_r}{\OPT}\ge\frac{\gamma^2}
{D_\gamma+\sqrt{15\gamma(1-\gamma)(r-1)}}$.
\end{proof}

\subsection{Proof of Theorem~\ref{thm:uniform-lower}}
\label{app:proof-uniform-lower}
\addappendixsubsection{Proof of Theorem~\ref{thm:uniform-lower}}
\begin{proof}
The first three terms in Eq.\eqref{eq:uniform-lower} come from Eq.\eqref{add_1}.  We next prove the final term in Eq.\eqref{eq:uniform-lower}. 
If $0<\gamma\le1/2$, then $2(1-\gamma)\ge1$, and hence
\[
B_\gamma+\sqrt{\gamma(r-1)}\le B_\gamma+\sqrt{2\gamma(1-\gamma)(r-1)} \le\sqrt2\,(B_\gamma+\sqrt{\gamma(1-\gamma)(r-1)}).
\]Thus, we can show that, when $0<\gamma\le1/2$,
\[R_{\gr}(\gamma,r)\ge\max\left\{\frac{\gamma}{r},
\frac{\gamma^2}{B_\gamma+\sqrt{\gamma(r-1)}}\right\}\ge\frac1{\sqrt2}\frac{\gamma^2}
{\left(B_\gamma+\sqrt{\gamma(1-\gamma)(r-1)}\right)}
\ge\frac{\gamma^2}
{14\left(B_\gamma+\sqrt{\gamma(1-\gamma)(r-1)}\right)},
\] where the first inequality follows from Eq~\eqref{eq:main-lower}.

It remains to consider $1/2\le\gamma<1$.  Write
$\delta:=1-\gamma$.  The universal rank-dependent term satisfies
\[
\frac{\gamma^2}{B_\gamma+\sqrt{\gamma(r-1)}}
\ge\frac{\gamma^2\sqrt{1-\gamma}}{B_\gamma+\sqrt{\gamma(1-\gamma)(r-1)}}.
\]
Thus, when $\sqrt{1-\gamma}\ge1/14$, we also can get 
\begin{equation*}
  R_{\gr}(\gamma,r)\ge\max\left\{\frac{\gamma}{r},
\frac{\gamma^2}{B_\gamma+\sqrt{\gamma(r-1)}}\right\}\ge\frac{\gamma^2}
{14\left(B_\gamma+\sqrt{\gamma(1-\gamma)(r-1)}\right)}.
\end{equation*}

Now, we suppose that $\sqrt{1-\gamma}<1/14$. Recall the result of Theorem~\ref{thm:near-lower}: when $1/2\le\gamma<1$, we can show 
\begin{equation*}
R_{\gr}(\gamma,r)
\ge
\max\left\{\frac{\gamma}{r},
\frac{\gamma^2}
{D_\gamma+\sqrt{15\gamma(1-\gamma)(r-1)}}\right\},
\end{equation*} where $D_\gamma\triangleq60-15\gamma+\gamma^2$. Next, we compare $D_\gamma+\sqrt{15\gamma(1-\gamma)(r-1)}$ with $14\big(B_\gamma+\sqrt{\gamma(1-\gamma)(r-1)}\big)$. Specifically, we can show that
\begin{equation*}
    \begin{aligned}
&D_\gamma+\sqrt{15\gamma(1-\gamma)(r-1)}-14\big(B_\gamma+\sqrt{\gamma(1-\gamma)(r-1)}\big)\\&=18+27(1-\gamma)-13(1-\gamma)^2-(14-\sqrt{15})\sqrt{\gamma(1-\gamma)(r-1)}.
    \end{aligned}
\end{equation*}If setting $u\triangleq\sqrt{\gamma(1-\gamma)(r-1)}$ and $u\ge 9/5$, we can have 
\begin{align*}
18+27(1-\gamma)-13(1-\gamma)^2<18+27(1-\gamma)<18+\frac{27}{196}
<(14-\sqrt{15})\frac{9}{5}\le (14-\sqrt{15})u,
\end{align*}where the second inequality follows from $\sqrt{1-\gamma}<1/14$ and  the third inequality uses
$\sqrt{15}<31/8$ and $(14-31/8)(9/5)=729/40=18+\frac{9}{40}\ge18+\frac{27}{196}$. Consequently,
\[D_\gamma+\sqrt{15\gamma(1-\gamma)(r-1)}\le 14\big(B_\gamma+\sqrt{\gamma(1-\gamma)(r-1)}\big)\] such that
\[
R_{\gr}(\gamma,r)
\ge
\max\left\{\frac{\gamma}{r},
\frac{\gamma^2}
{D_\gamma+\sqrt{15\gamma(1-\gamma)(r-1)}}\right\}\ge\frac{\gamma^2}
{14\left(B_\gamma+\sqrt{\gamma(1-\gamma)(r-1)}\right)}.
\]

Finally, we suppose that $u<9/5$ and set $x\triangleq\gamma^{r+1}$. Then, we can show that
\begin{equation*}
    \begin{aligned}
     \frac{\gamma^{r+1}}{2-\gamma+\gamma^{r+1}}
\bigg/\frac{\gamma^2}{B_\gamma+u}=\frac{\gamma^{r+1}(B_\gamma+u)}
{\gamma^2(2-\gamma+\gamma^{r+1})}= \frac{x(B_\gamma+u)}
{\gamma^2(2-\gamma+x)}.
    \end{aligned}
\end{equation*}Moreover,  $\sqrt{1-\gamma}<1/14$ implies $\gamma>195/196$. With this result, we also can have that
\[
B_\gamma\triangleq2-\gamma+\gamma^{2}\ge b_0\triangleq 2-\frac1{196}+\frac1{196^2},
\qquad 2-\gamma\le c_0\triangleq\frac{197}{196}.
\]
Therefore, we can get
\begin{align}
\frac{\gamma^{r+1}}{2-\gamma+\gamma^{r+1}}
\bigg/
\frac{\gamma^2}{B_\gamma+u}=\frac{x(B_\gamma+u)}
{\gamma^2(2-\gamma+x)}\ge\frac{x(B_\gamma+u)}
{2-\gamma+x}\ge\frac{x(b_0+u)}{c_0+x}.
\label{eq:endpoint-profile-comparison}
\end{align}
The inequality
$\log\gamma\ge-(1-\gamma)/\gamma$ gives
\begin{equation}\label{eq:endpoint-x-lower}
-\log x=-(r+1)\log\gamma
\le\frac{(1-\gamma)(r+1)}{\gamma}=\frac{\gamma(1-\gamma)(r-1)}{\gamma^{2}}+\frac{2(1-\gamma)}{\gamma}
=\frac{u^2}{\gamma^2}+\frac{2(1-\gamma)}{\gamma}.
\end{equation}
Since $u<9/5$, Eq.\eqref{eq:endpoint-x-lower} yields $-\log x
<\frac{81}{25}\left(\frac{196}{195}\right)^2
+\frac2{195}<\log27$. so $x>1/27$. From the monotonicity of
$x/(c_0+x)$ w.r.t. $x$, we can have 
\[
\frac{\gamma^{r+1}}{2-\gamma+\gamma^{r+1}}
\bigg/
\frac{\gamma^2}{B_\gamma+u}\ge\frac{x(b_0+u)}{c_0+x}
>\frac{b_0+u}{27c_0+1}.
\]If $u\ge1/64$, we can verify 
\begin{align*}
\frac{b_0+u}{27c_0+1}\ge \frac{b_0+1/64}{27c_0+1}=
\frac1{14}+\frac{109}{4323760}>
\frac1{14}
\end{align*}such that
\[
R_{\gr}(\gamma,r)
\ge\frac{\gamma^{r+1}}{2-\gamma+\gamma^{r+1}}\ge\frac{\gamma^2}
{14\left(B_\gamma+\sqrt{\gamma(1-\gamma)(r-1)}\right)}.
\]

If $u<1/64$, then Eq.\eqref{eq:endpoint-x-lower} gives $-\log x
<\frac1{4096}\left(\frac{196}{195}\right)^2
+\frac2{195}<\log2$. So $x>1/2$.  Hence
\[
\frac{x(b_0+u)}{c_0+x}
>\frac{b_0}{2c_0+1}
\approx0.662721
>\frac1{14},
\] which indicates that
\[
R_{\gr}(\gamma,r)
\ge\frac{\gamma^{r+1}}{2-\gamma+\gamma^{r+1}}\ge\frac{\gamma^2}
{14\left(B_\gamma+\sqrt{\gamma(1-\gamma)(r-1)}\right)}.
\]As a result, for every $0<\gamma\le1$ and every integer $r\ge2$, we have  
\begin{equation*}
    R_{\gr}(\gamma,r)\ge 
\frac{\gamma^2}
{14\left(B_\gamma+\sqrt{\gamma(1-\gamma)(r-1)}\right)},
\end{equation*}which proves Eq.\eqref{eq:uniform-lower} in Theorem~\ref{thm:uniform-lower}.
\end{proof}

\end{document}